\documentclass[12pt]{amsart}
\pdfoutput=1
\usepackage{amsmath,amsfonts,amssymb,amsthm,txfonts,pxfonts,amscd}
\usepackage{graphicx}
\usepackage[all]{xy}
\usepackage{subfigure}
\usepackage{geometry}
\usepackage{color}
\usepackage[T1]{fontenc}
\usepackage{a4wide}
\usepackage{bm}
\usepackage[mathcal]{euscript}

\graphicspath{{figs/}}

\newcommand{\mathnotation}[2]{\newcommand{#1}{\ensuremath{#2}}}

\renewcommand{\l}{\left}			
\renewcommand{\r}{\right}			
\mathnotation{\ldef}{\mathrel{\raisebox{.069ex}{:}\!\!=}}
\mathnotation{\rdef}{\mathrel{=\!\!\raisebox{.069ex}{:}}}
\mathnotation{\ho}{\varphi}
\mathnotation{\hoii}{\psi}
\mathnotation{\surf}{S}		
\mathnotation{\fol}{\mathcal{F}}
\mathnotation{\folu}{\fol^{\mathrm{u}}}
\mathnotation{\fols}{\fol^{\mathrm{s}}}
\mathnotation{\id}{Id}
\mathnotation{\tmeas}{\mu}
\mathnotation{\tmeasu}{\tmeas^{\mathrm{u}}}
\mathnotation{\tmeass}{\tmeas^{\mathrm{s}}}
\mathnotation{\Nsing}{N}
\mathnotation{\prongs}{\#\text{prongs}}

\mathnotation{\aaa}{a}				
\newcommand{\infloop}[1]{}			
\mathnotation{\fp}{p}
\mathnotation{\fq}{p'}
\mathnotation{\fpp}{S}
\mathnotation{\fpq}{Q}
\mathnotation{\xx}{X}				%

\newcommand{\R}{\mathbb{R}}

\newcommand{\Z}{\mathbb{Z}}

\newcommand{\Tr}{Tr}

\newcommand\D{d\hspace{-0.5pt}}
\newcommand\pA{\phi}
\newcommand\homeo{f}
\newcommand\homsing{p}
\newcommand{\DT}[1]{T_{#1}^{}}
\newcommand{\DTi}[1]{T_{#1}^{-1}}

\newtheorem{Theorem}{Theorem}[section]

\newtheorem{Question}{Question}[section]
\newtheorem{Corollary}[Theorem]{Corollary}
\newtheorem{Proposition}[Theorem]{Proposition}
\newtheorem{Lemma}[Theorem]{Lemma}
\newtheorem{Remark}{Remark}[section]
\newtheorem{Construction}{Construction}[section]

\newtheorem*{NoNumberTheorem}{Theorem}

\newtheorem*{Acknowledgments}{Acknowledgments}
\newtheorem*{NoNumberRemark}{Remark}
\newtheorem*{Convention}{Convention}

\theoremstyle{definition}
\newtheorem{Step}{Step}

\begin{document}

\title[Systole in small genus]
{On the minimum dilatation of pseudo-Anosov homeomorphisms on surfaces 
of small genus}
\author{Erwan Lanneau, Jean-Luc Thiffeault}

\address{
Centre de Physique Th\'eorique (CPT), UMR CNRS 6207 \newline
Universit\'e du Sud Toulon-Var and \newline
F\'ed\'eration de Recherches des Unit\'es de
Math\'ematiques de Marseille \newline
Luminy, Case 907, F-13288 Marseille Cedex 9, France
}
\email{erwan.lanneau@cpt.univ-mrs.fr}

\address{
Department of Mathematics \newline Van Vleck Hall, 480 Lincoln Drive \newline
 University of Wisconsin -- Madison, WI 53706, USA}

\email{jeanluc@math.wisc.edu}
\date{\today}


\begin{abstract}
  We find the minimum dilatation of pseudo-Anosov homeomorphisms that
  stabilize an orientable foliation on surfaces of genus three, four,
  or five, and provide a lower bound for genus six to eight. Our
  technique also simplifies Cho and Ham's proof of the least
  dilatation of pseudo-Anosov homeomorphisms on a genus two surface.
  For genus $g=2$ to $5$, the minimum dilatation is the smallest
  Salem number for polynomials of degree $2g$.
\end{abstract}

\maketitle

\renewcommand{\abstractname}{R\'esum\'e}

\begin{abstract}
  Nous calculons la plus petite dilatation d'un hom\'eomorphisme de
  type pseudo-Anosov laissant invariant un feuilletage mesur\'e
  orientable sur une surface de genre $g$ pour $g=3,4,5$. Nous donnons
  aussi une borne inf\'erieure pour les genres $6,7$ et $8$. Nos
  techniques simplifient la preuve de Cho et Ham sur le calcul de la
  plus petite dilatation d'un hom\'eomorphisme de type pseudo-Anosov
  sur une surface de genre $2$. Pour $g=2$ \`a $5$, la plus petite
  dilatation est le plus petit nombre de Salem pour les polynomes \`a
  degr\'e fix\'e $2g$.
\end{abstract}

\section{Introduction}

This paper concerns homeomorphisms of a compact oriented surface~$M$
to itself.  There are natural equivalence classes of such
homeomorphisms under isotopy, called isotopy classes or mapping
classes.  An irreducible mapping class is such that no power of its
members preserves a nontrivial subsurface.  By the Thurston--Nielsen
classification~\cite{Thurston1988}, irreducible mapping classes are
either finite-order or are of a type called pseudo-Anosov. The class
of pseudo-Anosov homeomorphisms is by far the richest.  One can think
of such a homeomorphism $\pA$ as an Anosov (or hyperbolic)
homeomorphism on $M \backslash \left\{ \textrm{singularities}
\right\}$. In particular, as for standard Anosov on the two
dimensional torus, there exists a local Euclidean structure (with
singularities) and two linear foliations ($\mathcal F^{s}$ and
$\mathcal F^{u}$, called stable and unstable) such that $\pA$ expands
the leaves of one foliation with a coefficient $\lambda$, and shrinks
those of the other foliation with the same coefficient. The
number~$\lambda$ is a topological invariant called the
\emph{dilatation} of $\pA$; the number $\log(\lambda)$ is the
\emph{topological entropy} of $\pA$.

Thurston proved that $\lambda +\lambda^{-1}$ is an algebraic integer
(in fact, it is a Perron number) over $\mathbb Q$ of degree bounded by
$4g-3$.  In particular Newton's formulas imply that for each
  $g\ge 2$ the set of dilatations bounded from above by a constant is
  finite. Hence the minimum value $\delta_g$ of the dilatation of
  pseudo-Anosov homeomorphisms on $M$ is well
  defined~\cite{Arnoux:Yoccoz, Ivanov}.  It can be shown that the
logarithm of~$\delta_g$ is the length of the shortest geodesics on the
moduli space of complex curves of genus $g$, $\mathcal M_g$ (for the
Teichm\"uller metric).

Two natural questions arise. The first is how to compute $\delta_g$
explicitly for small $g\geq 2$. The second question asks if there is a
unique (up to conjugacy) pseudo-Anosov homeomorphism with minimum
dilatation in the modular group Mod$(g)$.  It is well known that
$\delta_{1} = \frac1{2}(3+\sqrt{5})$ and this dilatation is uniquely
realized by the conjugacy class in Mod$(1)=\textrm{PSL}_2(\Z)$ of the
matrix $\left(\begin{smallmatrix}2&1\\1&1\end{smallmatrix}\right)$. In
principle these dilatations can be computed for any given $g$ using
train tracks. Of course actually carrying out this procedure, even for
small values of $g$, seems impractical.

We know very little about the value of the constants~$\delta_g$.
Using a computer and train tracks techniques for the punctured disc,
Cho and Ham~\cite{Ham2008} proved that $\delta_2$ is equal to the
largest root of the polynomial $X^4-X^3-X^2-X+1$, $\delta_{2} \simeq
1.72208$~\cite{Ham2008}.  One of the results of the present paper is
an independent and elementary proof of this fact.

One can also ask about the uniqueness (up to conjugacy) of
pseudo-Anosov homeomorphisms that realize $\delta_g$.  In genus $2$,
$\delta_2$ is not unique due to the existence of the hyperelliptic
involution and covering transformations (see Section~\ref{sec:genus2}
and Remark~\ref{rem:unicity} for a precise definition). But, up to
hyperelliptic involution and covering transformations, we prove the
uniqueness of the conjugacy class of pseudo-Anosov homeomorphisms that
realize $\delta_2$, in the mapping class group of genus $2$ surfaces,
Mod$(2)$ (see Theorem~\ref{thm:genre2}).

For $g > 1$ the estimate $2^{1/(12g-12)} \le \delta_g \le
(2+\sqrt{3})^{1/g}$ holds~\cite{Penner1991,Hironaka2006}.  We will
denote by $\delta^{+}_g$ the minimum value of the dilatation of
pseudo-Anosov homeomorphisms on a genus $g$ surface with
\emph{orientable} invariant foliations.  We shall prove:

\begin{Theorem}
\label{thm:genre2}
The minimum dilatation of a pseudo-Anosov homeomorphism 
on a genus two surface is equal to the largest root of the polynomial 
$X^4-X^3-X^2-X+1$,
\[
\delta_2 = \delta^{+}_2 = \tfrac{1}{4}+\tfrac{\sqrt{13}}{4}+
\tfrac{1}{2}\sqrt{\tfrac{\sqrt{13}}{2}-\tfrac{1}{2}} \simeq 1.72208.
\]
Moreover there exists a unique (up to conjugacy, hyperelliptic
involution, and covering transformations) pseudo-Anosov homeomorphism
on a genus two surface with dilatation $\delta_2$.
\end{Theorem}

\begin{NoNumberRemark}
This answers Problem~$7.3$ and Question~$7.4$ of Farb~\cite{Farb} in genus two.
\end{NoNumberRemark}

\begin{Theorem}
\label{thm:1}
The minimum value of the dilatation of pseudo-Anosov homeomorphisms on
a genus $g$ surface, $3\le g\le 5$, with orientable invariant
foliations is equal to the largest root of the polynomials in
Table~\ref{tab:thm:systoles}.
\end{Theorem}
\begin{table}[htbp]
\begin{tabular}{lll}
\hline
$g$ & polynomial  & $\delta^{+}_g\simeq$ \\
\hline
3 & $X^6 - X^4 - X^3 - X^2 + 1$ &  1.40127 \\
4 & $X^8-X^5-X^4-X^3+1$ & 1.28064 \\
5 & $X^{10} + X^9 - X^7 - X^6 - X^5 - X^4 - X^3 + X + 1$ & 1.17628 \\
\hline
\end{tabular}
\medskip
\caption{}
\label{tab:thm:systoles}
\end{table}
All of the minimum dilatations for~$2\le g \le 5$ are Salem
numbers~\cite{PisotSalem}.  In fact, their polynomials have the smallest
Mahler measure over polynomials of their degree~\cite{Boyd1980}.
For~$g=5$, the dilatation is realized by the pseudo-Anosov
homeomorphism described by Leininger~\cite{Leininger2004} as a
composition of Dehn twists about two multicurves.  Its characteristic
polynomial is the irreducible one having Lehmer's number as a root:
this is the smallest known Salem number.  The polynomial has the
smallest known Mahler measure over all integral polynomials.

For~$g=3$ and~$4$, we have constructed explicit examples. We present
two independent constructions in this paper: The first is given in
term of Dehn twists on a surface; The second involves the
\emph{Rauzy--Veech} construction (see
Appendix~\ref{appendix:construction:pA}).

\begin{Theorem}
\label{thm:lower}
The minimum value of the dilatation of pseudo-Anosov homeomorphisms
on a genus $g$ surface, $6\le g\le 8$, with orientable invariant
foliations is not less than the largest root of the polynomials in
Table~\ref{tab:thm:lower}. \\
In particular $\delta_6^{+} \geq \delta_5^{+}$.
\end{Theorem}
\begin{table}[htbp]
\begin{tabular}{lll}
\hline
$g$ & polynomial  & $\delta^{+}_g\gtrsim$  \\
\hline
6 & $X^{12} - X^7 - X^6 - X^5 + 1$ &  1.17628   \\
7 & $X^{14} + X^{13} - X^9 - X^8 - X^7 - X^6 - X^5 + X + 1$ & 1.11548 \\
8 & $X^{16} - X^9 - X^8 - X^7 + 1$ & 1.12876 \\
\hline
\end{tabular}
\medskip
\caption{}
\label{tab:thm:lower}
\end{table}

\begin{Remark}
Genus 6 is the first instance of a nondecreasing dilatation
compared to the previous genus.  This partially answers Question~$7.2$ of
Farb~\cite{Farb} in the orientable case.

We have also found an example of a pseudo-Anosov homeomorphism on a
genus $3$ surface that stabilizes a non-orientable measured foliation,
with dilatation $\delta_3^+$.  There is also evidence that~$\delta_5 <
\delta_5^+$~\cite{Aaber2010} (Section~\ref{sec:proofg5}).  In
addition, Aaber \& Dunfield~\cite{Aaber2010} and Kin \&
Takasawa~\cite{Kin2010} have found a pseudo-Anosov homeomorphism
realizing~$\delta_7^+$, and Hironaka~\cite{Hironaka2009} has done the
same for~$\delta^+_8$.  Hence, all the lower bounds in
Table~\ref{tab:thm:lower} except for genus~$6$ are known to be
realized by a pseudo-Anosov homeomorphism.

\end{Remark}

\begin{Remark}
  Our techniques also provide a way to investigate least dilatations
  of punctured discs. This will appear in the forthcoming
  paper~\cite{Lanneau:Thiffeault}.  Note that, for genus 3 to 8, none
  of the minimum dilatations realizing the bound can come from the
  lift of a pseudo-Anosov on a punctured disk (or any other
  lower-genus surface). Indeed, if the pseudo-Anosov comes from a
  lift, then composing this pseudo-Anosov with the hyperelliptic
  involution, one gets two pseudo-Anosov homeomorphisms, one with
  positive root when acting on homology, and one with negative root.
  Since the polynomials we find have only one sign of the dominant
  root when acting on homology, a lift is always ruled out.  This is
  in contrast to the Hironaka \& Kin~\cite{Hironaka2006} examples,
  which come from punctured disks.
\end{Remark}

\begin{Acknowledgments}
  The authors thank Christopher Leininger, Fr\'ed\'eric Le Roux,
  J\'er\^ome Los, Sarah Matz, and Rupert Venzke for helpful
  conversations, and are grateful to Matthew D.~Finn for help in
  finding pseudo-Anosov homeomorphisms in terms of Dehn twists.  J-LT
  thanks the Centre de Physique Th\'eorique de Marseille, where this
  work began, for its hospitality.  J-LT was also supported
  by the Division of Mathematical Sciences of the US National Science
  Foundation, under grant DMS-0806821.
\end{Acknowledgments}

\section{Background and tools}
\label{sec:background}

In this section we recall some general properties of dilatations and
pseudo-Anosov homeomorphisms, namely algebraic and spectral radius
properties. We also summarizes basic tools for proving our results
(for example see~\cite{Thurston1988,FLP,MT,Mc}).

To guide the reader, we will first outline the general method used to
find the least dilatation~$\delta_g^+$:\medskip

\emph{
Summary: to find the least dilatation~$\delta_g^+$ on a surface~$M$ of
genus~$g$.
\begin{enumerate}
\item Start with a known pseudo-Anosov homeomorphism on~$M$, with
  dilatation~$\alpha$, that stabilizes orientable foliations (we use
  the family in~\cite{Hironaka2006}).
\item Enumerate all reciprocal polynomials with Perron root less
  that~$\alpha$ (see Section~\ref{sec:algprop} for definitions, and
  Appendix~\ref{sec:smallP} for an explicit algorithm).  For genus~$g>2$,
  this requires a computer, but is a standard calculation.
\item Of these polynomials, eliminate the ones that are incompatible
  with the Lefschetz theorem (see Section~\ref{sec:Lef}).  The
  remaining polynomial with the smallest root gives a lower bound on
  the least dilatation~$\delta_g^+$.  For genus~$g>4$, this step
  requires a computer.
\item If possible, construct an explicit pseudo-Anosov homeomorphism
  on~$M$ having the lower bound in the previous step as a dilatation.
  We do this by either exhibiting a sequence of Dehn twists, or by the
  Rauzy--Veech construction (see
  Appendix~\ref{appendix:construction:pA}).  This confirms that we
  have found~$\delta_g^+$.
\end{enumerate}
}

\subsection{Affine structures and affine homeomorphisms}

To each pseudo-Anosov homeomorphism $\pA$ one can associate an affine
structure on $M$ for which $\pA$ is affine.

\subsubsection{Affine structures}
A surface of genus $g\geq 1$ is called a flat surface if it can be
obtained by edge-to-edge gluing of polygons in the plane using
translations or translations composed with $-\id$. We will call such a
surface $(M,q)$ where $q$ is the form $\D z^2$ defined locally.  The
metric on $M$ has zero curvature except at the zeroes of $q$ where the
metric has \emph{conical singularities} of angle $(k+2)\pi$ (with
$k\geq -1$). The integer $k$ is called the \emph{degree} of the zero
of $q$. A point that is not singular is {\it regular}. We will use the
convention that a singular point of degree $0$ is regular. A measured
foliation $M$ is a linear flow on this flat surface $M$ for an affine
structure.

The Gauss--Bonnet formula applied to the singularities reads $\sum_i k_i = 4g-4$.
We will call the integer vector (or simply the stratum) 
$(k_1,\dots,k_n)$ with $k_i\geq-1$ the \emph{singularity data} of the measured foliation. 

If one restricts gluing to translations only then the surface is
called a translation surface; otherwise it is called a
half-translation surface.  For a translation surface the degree of all
singularities is even; the converse is false in general.

There is a standard construction, {\it the orientating cover}, that
produce a translation surface from a half-translation surface.

\begin{Construction}
\label{rk:orientating:cover}
Let $N$ be a half-translation surface with singularity data
$(k_1,\dots,k_n)$. Then there exists a translation surface $M$ and a
double branched cover $\pi : M \rightarrow N$, branched precisely over
the singular points of odd degree. In addition $\pi$ is the minimal
double branched cover in this class.
\end{Construction}

\subsubsection{Affine homeomorphisms}
A homeomorphism $f$ is affine with respect to $(M,q)$ if $f$ permutes
the singularities, $f$ is a diffeomorphism on the complement of the
singularities, and the derivative map $D f$ of $f$ is a constant
matrix in $\textrm{PSL}_2(\R)$.

There is a standard classification of the elements of
$\textrm{PSL}_2(\R)$ into three types: elliptic, parabolic and
hyperbolic. This induces a classification of affine homeomorphisms. An
affine homeomorphism is parabolic, elliptic, or pseudo-Anosov,
respectively, if $|\Tr(D f)| = 2$, $\Tr(D f)| < 2$, or $|\Tr(D f)| >
2$, respectively (where~$\Tr$ is the trace).

\subsubsection{Pseudo-Anosov homeomorphisms}

Since we are interested in pseudo-Anosov homeomorphisms we will assume
that $|\Tr(D f)| > 2$.  Then there exists an eigenvalue $\lambda$ of
$D f$ such that $|\lambda|>1$ and $\Tr(D f) = \lambda + \lambda^{-1}$.
The two eigenvectors associated to $\lambda$ and $\lambda^{-1}$
determine two directions on the flat surface $M$, invariant by
$\pA$. Of course $\pA$ expends leaves of the stable foliation by the
factor $|\lambda|$ and shrinks leaves of the unstable foliation by the
same factor. We can assume that these directions are horizontal and
vertical. In these coordinates $(M,q)$, the pair of associated
measured foliations (stable and unstable) of $\pA$ are given by the
horizontal and vertical measured foliations Im$(q)$ and Re$(q)$ and
the derivative of $\pA$ is the matrix $A=\left( \begin{smallmatrix}
    \pm \lambda^{-1} & 0 \\ 0 & \pm \lambda \end{smallmatrix}
\right)$.  By construction the dilatation $\lambda(\pA)$ of $\pA$
equals $|\lambda|$. The singularity data of a pseudo-Anosov $\pA$ is
the singularity data of its invariant measured foliation.

The group $\textrm{PSL}_2(\mathbb R)$ naturally acts on the set of
flat surfaces.  With above notations the matrix $A$ fixes the surface
$(M,q)$, that is, $(M,q)$ can be obtained from $A \cdot (M,q)$ by
``cutting'' and ``gluing'' (i.e. the two surfaces represent the same
point in the moduli space).  The converse is true: if $A$ stabilizes a
flat surface $(M,q)$, then there exists an affine diffeomorphism $f:M
\rightarrow M$ such that $D f=A$. \medskip

Masur and Smillie~\cite{Masur:Smillie} proved the following result:

\begin{Theorem}[Masur,~Smillie]
  For each integer partition $(k_1,\dots,k_n)$ of $4g-4$ with
  $k_i\geq0$ even, there is a pseudo-Anosov homeomorphism $\pA$ with
  singularity data $(k_1,\dots,k_n)$ that fixes an orientable measured
  foliation.
  For each integer partition $(k_1,\dots,k_n)$ of $4g-4$ with
  $k_i\geq-1$, there is a pseudo-Anosov homeomorphism $\pA$ with
  singularity data $(k_1,\dots,k_n)$ that fixes a non-orientable
  measured foliation, with the following exceptions:
\[
(1,-1), \ (1,3),\textrm{ and } (4).
\]
\end{Theorem}

\begin{Convention}
  For the remainder of this paper, unless explicitly stated (in
  particular in Section~\ref{sec:genus2}), we shall assume that
  pseudo-Anosov homeomorphisms preserve \emph{orientable} measured
  foliations.
\end{Convention}

For instance, if $g=3$ and $\pA$ preserves an orientable measured
foliation, then there are $5$ possible strata for the singularity data
of $\pA$:
\[
(8), \ (2,6), \ (4,4), \ (2,2,4),\textrm{ and } (2,2,2,2).
\]

\subsection{Algebraic properties of dilatations}
\label{sec:algprop}

The next theorem follows from basic results in the theory of
pseudo-Anosov homeomorphisms (see for example~\cite{Thurston1988}).

\begin{Theorem}[Thurston]
\label{thm:thurston}
Let $\pA$ be a pseudo-Anosov homeomorphism on a genus $g$ surface that
leaves invariant an orientable measured foliation. Then
\begin{enumerate}

\item The linear map $\pA_\ast$ defined on $H_1(M,\mathbb R)$ has a
  simple eigenvalue $\rho(\pA_\ast) \in \mathbb R$ such that
  $|\rho(\pA_\ast)| > |x| $ for all other eigenvalues $x$;

\item $\pA$ is affine, for the affine structure determined by the
  measured foliations, and the eigenvalues of the derivative $D\pA$
  are $\rho(\pA_\ast)^{\pm 1}$;

\item $|\rho(\pA_\ast)| > 1$ is the dilatation $\lambda$ of $\pA$. 

\end{enumerate}

\end{Theorem}

A \emph{Perron root} is an algebraic integer $\lambda \geq 1$ all whose
other conjugates satisfy $|\lambda'| < \lambda$. Observe that these
are exactly the numbers that arise as the leading eigenvalues of
Perron--Frobenius matrices.  Since $\pA_\ast$ preserves a symplectic
form, the characteristic polynomial $\chi_{\pA_\ast}$ is a reciprocal
degree $2g$ polynomial.

\begin{Remark}
\label{rk:perron}
The dilatation of a pseudo-Anosov homeomorphism $\pA$ is the Perron
root of a reciprocal degree $2g$ polynomial, namely
$\chi_{\pA_\ast}(X)$ if $\rho(\pA_\ast)>0$ and $\chi_{\pA_\ast}(-X)$
otherwise.
\end{Remark}

There is a converse to Theorem~\ref{thm:thurston}, but the proof does
not seem as well-known, so we include a proof here
(see~\cite{Band:Boyland} Lemma~$4.3$).

\begin{Theorem}
\label{theo:orientation:converse}
Let $\pA$ be a pseudo-Anosov homeomorphism on a surface $M$ with
dilatation $\lambda$. Then the following are equivalent:
\begin{enumerate}

\item $\lambda$ is an eigenvalue of the linear map $\pA_\ast$ defined
  on $H_1(M,\mathbb R)$.

\item The invariant measured foliations of $\pA$ are orientable.

\end{enumerate}
\end{Theorem}

\begin{proof}
  Suppose the stable measured foliation on $(M,q)$ is
  non-orientable. There exists a double branched cover $\pi: N
  \rightarrow M$ which orients the foliation (we denote by $\tau$ the
  involution of the covering).  Let $[w]$ be an eigenvector of
  $\pA_\ast$ in $H^1(M,\R)$ with eigenvalue $\lambda$. The vector
  $[w]$ pulls back to an eigenvector $[w']$ of the adjoint $\phi^\ast$
  in $H^1(N,\R)$ for the eigenvalue $\lambda$.

  The stable foliation on $N$ now also defines a cohomology class
  $[Re(\omega)]$ where $\omega^2=\pi^\ast q$.  By construction
  $[Re(\omega)]$ is an eigenvector for the eigenvalue $\lambda$. By
  Theorem~\ref{thm:thurston} $\lambda$ is simple so that the two
  classes $[Re(\omega)]$ and $[w']$ must be linearly dependent.  But
  since $[w']$ is invariant by the deck transformation $\tau$, while
  $[Re(\omega)]$ is sent to $-[Re(\omega)]$ by $\tau$, we get a
  contradiction.
\end{proof}

Combining this theorem with two classical results of
Casson--Bleiler~\cite{Casson-Bleiler} and Thurston~\cite{FLP} we get
\begin{Theorem}
\label{thm:construction:pA}
Let $f$ be a homeomorphism on a surface $M$ and let $P(X)$ be the
characteristic polynomial of the linear map $f_\ast$ defined on
$H_1(M,\mathbb R)$. Then one has
\begin{enumerate}

\item If $P(X)$ is irreducible over $\Z$, has no roots of unity as
  zeroes, and is not a polynomial in $X^k$ for $k>1$, then $f$ is
  isotopic to a pseudo-Anosov homeomorphism $\pA$;
\label{thm:construction:pA:1}
\item In addition, if the maximal eigenvalue (in absolute value) of
  the action of $f$ on the fundamental group is $\lambda>1$, then the
  dilatation of $\pA$ is $\lambda$;
\item In addition, if $\lambda$ is the Perron root of $P(X)$, then
  $\pA$ leaves invariant orientable measured foliations.
\end{enumerate}
\end{Theorem}

\begin{proof}
  The first point asserts that $f$ is isotopic to a pseudo-Anosov
  homeomorphism $\pA$~\cite[Lemma~$5.1$]{Casson-Bleiler}. The second
  point asserts that $\pA$ has dilatation $\lambda$
  \cite[Expos\'e~$10$]{FLP}.  Finally by the previous theorem, the
  last assumption implies that the invariant measured foliations of
  $\pA$ are orientable.
\end{proof}

We will need a more precise statement. The following has been remarked
by Bestvina:

\begin{Proposition}
\label{rem:Bestvina}
The statement ``$P$ is irreducible over $\Z$'' in part
(\ref{thm:construction:pA:1}) of Theorem~\ref{thm:construction:pA} can
be replaced by ``$P$ is symplectically irreducible over $\Z$'',
meaning that $P$ is not the product of two nontrivial reciprocal
polynomials.
\end{Proposition}

\subsection{Pseudo-Anosov homeomorphisms and the Lefschetz theorem}
\label{sec:Lef}

In this section, we recall the well-known Lefschetz theorem for
\emph{homeomorphisms} on compact surfaces (see for
example~\cite{Brown}).  If $\homsing$ is a fixed point of a
homeomorphism $\homeo$, we define the index of $\homeo$ at $\homsing$
to be the algebraic number $\textrm{Ind}(\homeo,\homsing)$ of turns of
the vector $(x,\homeo(x))$ when $x$ describes a small loop around
$\homsing$.

\begin{NoNumberTheorem}[Lefschetz theorem]

  Let $\homeo$ be a homeomorphism on a compact surface $M$.  Denote by
  $\Tr(\homeo_\ast)$ the trace of the linear map $\homeo_\ast$ defined
  on the first homology group $H_1(M,\mathbb R)$.  Then the Lefschetz
  number $L(\homeo)$ is $2-\Tr(\homeo_\ast)$. Moreover the following
  equality holds:
\[
L(\homeo) = \sum_{\homsing=f(\homsing)} \textrm{Ind}(\homeo,\homsing).
\]
\end{NoNumberTheorem}

For a pseudo-Anosov homeomorphism $\pA$, if $\Sigma \in M$ is a
singularity of the stable foliation of $\pA$ (of degree $2d$) then
there are $2(d+1)$ emanating rays. The orientation of the foliation
defines $d+1$ outgoing separatrices and $d+1$ ingoing separatrices.

\begin{Proposition}
\label{prop:index}

Let $\Sigma$ be a fixed singularity of $\pA$ of degree $2d$ and let
$\rho(\pA_\ast)$ be the leading eigenvalue of $\pA_\ast$.  Then
\begin{itemize}

\item If $\rho(\pA_\ast) < 0$ then $\pA$ exchanges the set of outgoing
separatrices and the set of ingoing separatrices. Moreover 
$\textrm{Ind}(\pA,\Sigma) = 1$.

\item If $\rho(\pA_\ast) > 0$ then either

\begin{itemize}

\item $\pA$ fixes each separatrix and $\textrm{Ind}(\pA,\Sigma) =
  1-2(d+1) <0$, or

\item $\pA$ permutes cyclically the outgoing separatrices (and 
ingoing separatrices) and $\textrm{Ind}(\pA,\Sigma) = 1$.

\end{itemize}
\end{itemize}

\end{Proposition}

\begin{figure}
\subfigure[]{
  \includegraphics[height=.2\textheight]{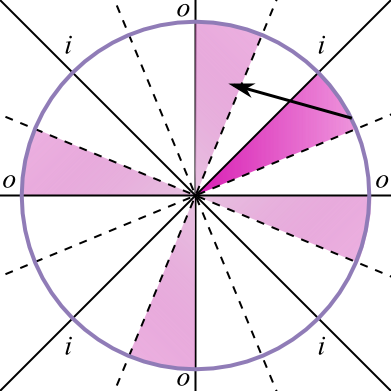}
}\hspace{3em}
\subfigure[]{
  \includegraphics[height=.2\textheight]{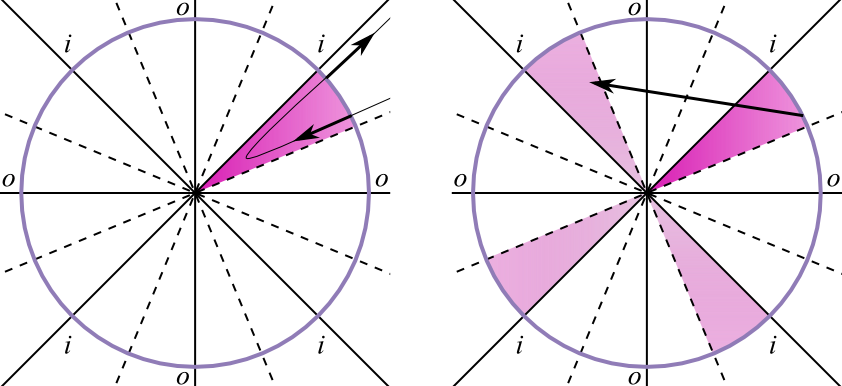}
}
\caption{Mapping of the~$4(d+1)$ hyperbolic sectors by~$\phi$ near a
  degree~$2d=6$ singularity: (a) $\rho(\pA_\ast)<0$: the sectors are
  permuted and the index is~$1$; (b) $\rho(\pA_\ast)>0$: the sectors
  can either be fixed (left, index~$1-2(d+1)=-7$) or permuted (right,
  index~$1$).  The index is defined as the number of turns of a vector
  joining $x$ to $\phi(x)$ as $x$ travels counterclockwise around a
  small circle.  The separatrices of the unstable foliation are
  alternately labeled ingoing ($i$) and outgoing ($o$).  The grey
  areas indicate a hyperbolic sector and its possible images for each
  case.}
\label{fig:index}
\end{figure}

\begin{proof}[Proof of Proposition~\ref{prop:index}]
  Obviously $\pA$ acts on the set of separatrices (namely the set of
  outgoing separatrices and ingoing separatrices). It is clear that
  $\rho(\pA_\ast) < 0$ if and only if $\pA$ exchanges these two
  sets. In that case, $\textrm{Ind}(\pA,\Sigma) = 1$ for any fixed
  point $\Sigma$, since the tip of the vector~$(x,\homeo(x))$ never
  crosses the hyperbolic sector containing~$x$ and is thus constrained
  to make a single turn counterclockwise.  (A hyperbolic sector is the
  region between adjacent ingoing and outgoing separatrices, see
  figure~\ref{fig:index}.)

  If $\rho(\pA_\ast) > 0$ then $\pA$ fixes globally the set of
  outgoing separatrices.  Let us assume that $\pA$ fixes an outgoing
  separatrix $\gamma^{u}$ of the unstable foliation $\mathcal
  F^{u}$. Let $\gamma_1^{s}$ and $\gamma_2^{s}$ be two adjacent
  incoming separatrices for the stable foliation $\mathcal F^{s}$ that
  define a sector containing~$\gamma^{u}$ and another (ingoing)
  separatrix of $\mathcal F^{u}$.  Since $\gamma^u$ is fixed by $\pA$,
  the sector determined by $\gamma_1^{s}$ and $\gamma_2^{s}$ is also
  fixed. $\pA$ preserves orientation so that $\gamma_1^{s}$ (and so
  $\gamma_2^{s}$) is fixed.  Hence, the other separatrix of $\mathcal
  F^{u}$ in the sector is fixed.  By induction, each separatrix of
  $\mathcal F^{u}$ is fixed.

  There are $4(d+1)$ hyperbolic sectors. For each sector, the vector
  $(x,h(x)$ describes an angle of $-\pi$ plus the sector angle,
  $\pi/2(d+1)$.  Thus the total angle is $4(d+1)(-\pi + \pi/2(d+1)) =
  2\pi (1-2(d+1))$.

  If $\pA$ has no fixed separatrices then clearly $\pA$ permutes the
  outgoing separatrices. In addition, $\pA$ is isotopic to a rotation,
  thus $\pA$ permutes cyclically the separatrices~\cite{Leroux2004}.
  In that case~$\textrm{Ind}(\pA,\Sigma) = 1$, for the same reason as
  the $\rho(\pA_\ast) < 0$ case above.
\end{proof}

We will use the following corollaries:

\begin{Corollary}[Lefschetz theorem revisited for pseudo-Anosov
  homeomorphisms]
\label{cor:LpA}

Let $Sing(\pA)$ be the set of fixed singularities of degree $d>0$ of
the pseudo-Anosov homeomorphism $\pA$. Let $Fix(\pA)$ be the set of
regular fixed points of $\pA$ (i.e. of degree $d=0$).

Then if $\rho(\pA_\ast) > 0$,
\[
2-\Tr(\pA_\ast) = \sum_{\Sigma \in Sing(\pA)}\textrm{Ind}(\pA,\Sigma) - \#Fix(\pA)
\]
where $\textrm{Ind}(\pA,\Sigma)=1$ or $1-2(d+1)$ and $2d$ is the degree of
$\Sigma$.
\medskip

If $\rho(\pA_\ast) < 0$,
\[
2-\Tr(\pA_\ast) = \#Sing(\pA) + \#Fix(\pA)\,.
\]

\end{Corollary}

\begin{Corollary}
\label{cor:trick}
Let $\Sigma$ be a fixed singularity of $\pA$ (of degree $2d$). Let us
assume that $\rho(\pA_\ast) > 0$ and
$\textrm{Ind}(\pA,\Sigma)=1$. Then
\[
\forall i=1,\dots d ,\ \textrm{Ind}(\pA^i,\Sigma)=1
\]
and 
\[
\textrm{Ind}(\pA^{d+1},\Sigma)=1-2(d+1).
\]
\end{Corollary}

We will use this corollary with $d=2$ and $d=4$ in the coming
sections, so we prove it only for those cases.

\begin{proof}[Proof of Corollary~\ref{cor:trick}]
  If $\Sigma$ is a singularity of degree $2$ ($d=1$) then there are
  $2$ outgoing separatrices.  $\textrm{Ind}(\pA,\Sigma)=1$ implies
  that $\pA$ permutes these two separatrices so that $\pA^2$ fixes
  them. Hence $\textrm{Ind}(\pA^2,\Sigma)=1-2(1+1)=-3$.

  If $\Sigma$ is a singularity of degree $4$ ($d=2$) then there are
  three outgoing separatrices.  $\textrm{Ind}(\pA,\Sigma)=1$ implies
  that $\pA$ permutes cyclically these three separatrices.  Hence
  $\textrm{Ind}(\pA^2,\Sigma)=1$ and
  $\textrm{Ind}(\pA^3,\Sigma)=1-2(2+1)=-5$.
\end{proof}

\section{Genus three: A proof of Theorem~\ref{thm:1} for~$g=3$}
\label{sec:proof:g3}

We write~$\rho(P)$ for the largest root (in absolute value) of a
polynomial~$P$; for the polynomials we consider this is always real
and with strictly larger absolute value than all the other roots,
though it could have either sign.  If~$\rho(P)>0$ then it is a Perron
root; otherwise~$\rho(P(-X))$ is a Perron root.

We find all reciprocal polynomials with a Perron root less than our
candidate and then we test whether a polynomial is compatible with a
given stratum.  This is straightforward: we simply try all possible
permutations of the singularities and separatrices, and calculate the
contribution to the Lefschetz numbers for each iterate of~$\phi$.
Then we see whether the deficit in the Lefschetz numbers can be
exactly compensated by regular periodic orbits.  If not, the
polynomial cannot correspond to a pseudo-Anosov homeomorphism on that
stratum.

We prove the theorems out of order since genus 3 is simplest. 
We know that $\delta^{+}_3 \leq \rho(X^3-X^2-1)\simeq 1.46557$ (for instance
see~\cite{Hironaka2006} or~\cite{Lanneau:Thiffeault}). We will construct 
a pseudo-Anosov homeomorphism with a smaller dilatation than~$1.46557$ 
and prove that this dilatation is actually the least dilatation.

Recall that $\delta^{+}_3$ is the Perron root of some reciprocal
polynomial $P$ of degree $6$ (see Remark~\ref{rk:perron}).  As
discussed in Appendix~\ref{sec:smallP}, it is not difficult to find
all reciprocal polynomials with a Perron root $\rho(P)$, $1 < \rho(P)
< \rho(X^3-X^2-1)$: there are only two, listed in
Table~\ref{tab:disc6} (see also Appendix~\ref{subsection:mahler} for an
alternate approach to this problem).
\begin{table}[htbp]
\begin{tabular}{lr}
\hline
polynomial  & Perron root  \\
\hline
$P_1 = (X^3-X-1)(X^3+X^2-1)$  &  1.32472  \\
$P_2 = X^6 - X^4 - X^3 - X^2 +1$ &  1.40127   \\
\hline
\end{tabular}
\medskip
\caption{List of all reciprocal monic degree $6$ polynomials $P$ with 
Perron root $1 < \rho(P) < \rho(X^3-X^2-1)\simeq 1.46557$.}
\label{tab:disc6}
\end{table}
Let us assume that $\delta^{+}_3 < \rho(X^3-X^2-1)$ and see if we get
a contradiction.  We let $\pA$ be a pseudo-Anosov homeomorphism with
$\lambda(\pA) = \delta^{+}_3$.  By the above discussion there are only
two possible candidates for a reciprocal annihilating polynomial $P$
of the dilatation of $\pA$, namely $\lambda(\pA) = \rho(P_i)$ for some
$i \in \{1,2\}$.  In the next subsection we shall prove that there are
no pseudo-Anosov homeomorphisms on a genus three surface (stabilizing
orientable foliations) with a dilatation $\rho(P_1)$.  We shall then
show that a pseudo-Anosov homeomorphism with dilatation~$\rho(P_2)$
exists on this surface.

\subsection{First polynomial: $\lambda(\pA) = \rho(P_1)$}

Let $\pA_\ast$ be the linear map defined on $H_1(X,\mathbb R)$ and let
$\chi_{\pA_\ast}$ be its characteristic polynomial. By
Theorem~\ref{thm:thurston} the leading eigenvalue $\rho(\pA_\ast)$ of
$\pA_\ast$ is $\pm \rho(P_1)$. The minimal polynomial of the
dilatation of $\pA$ is $X^3-X-1$; thus if $\rho(\pA_\ast) > 0$ then
$X^3-X-1$ divides $\chi_{\pA_\ast}$, otherwise $X^3-X+1$ divides
$\chi_{\pA_\ast}$.  Requiring the polynomial to be reciprocal leads to
$\chi_{\pA_\ast} = P_1$ for the the first case and $\chi_{\pA_\ast}
=P_1(-X)= (X^3-X+1)(X^3-X^2+1)$ for the second.

The trace of $\pA_\ast^n$ (and so the Lefschetz number of $\pA^n$) is
easy to compute in terms of its characteristic polynomial. Let us
analyze carefully the two cases depending on the sign of
$\rho(\pA_\ast)$.

\begin{enumerate}

\item If $\rho(\pA_\ast) < 0$ then $\chi_{\pA_\ast}(X) = P_1(-X) =
  (X^3-X+1)(X^3-X^2+1)$. Let $\psi=\phi^2$.  Observe that $\psi$ is a
  pseudo-Anosov homeomorphism and $\rho(\psi_\ast) > 0$ is a Perron
  root.  >From Newton's formulas (see Appendix~\ref{sec:smallP}), we
  have $\Tr(\phi_\ast)=-1$, $\Tr(\psi_\ast)=3$, $\Tr(\psi^2_\ast)=-1$,
  and~$\Tr(\psi^3_\ast)=3$, so that $L(\phi)=3$, $L(\psi)=-1$,
  $L(\psi^2)=3$, and~$L(\psi^3)=-1$.

  As we have seen in Section~\ref{sec:background}, there are $5$
  possible strata for the singularity data of $\pA$, and so for
  $\psi$, namely,
\[
(8), \ (2,6), \ (4,4), \ (2,2,4), \textrm{ and } (2,2,2,2).
\]
Since $L(\psi^2) = 3$ there are at least $3$ singularities (of
index~$+1$) fixed by $\psi^2$; thus we need only consider strata
$(2,2,4)$ and $(2,2,2,2)$.  (From Corollary~\ref{cor:LpA} regular
fixed points can only give negative index since~$\rho(\psi^2_*)>0$.)

For stratum~$(2,2,4)$, the single degree-4 singularity must be fixed,
and its three outgoing separatrices must be fixed by~$\psi^3$.  The
contribution to the index is then~$-5$, which
contradicts~$L(\psi^3)=-1$ since there is no way to make up the
deficit.

For stratum~$(2,2,2,2)$, since~$\psi^2$ fixes at least three
singularities they account for~$+3$ of the Lefschetz
number~$L(\psi^2)=3$.  But the fourth singularity must also be fixed
by~$\psi^2$, so it adds~$+1$ or~$-3$ to the Lefschetz number,
depending on the permutation of its two separatrices.  The only
compatible scenario is that it adds~$+1$, with the difference
accounted by a single regular fixed point that contributes~$-1$.
Since all four singularities are thus fixed by~$\psi^2=\phi^4$, this
means that their permutation~$\sigma\in S_4$ must
satisfy~$\sigma^4=\id$.  There are three cases: either the
singularities are all fixed by~$\phi$, they are permuted in groups of
two, or they are cyclically permuted.  For the first two cases, the
singularities are also fixed by~$\psi=\phi^2$, so by
Corollary~\ref{cor:trick} they cannot contribute positively
to~$\psi^2$, which they must as we saw above.  If the four
singularities are all cyclically permuted, then they contribute
nothing to~$L(\phi)=3$ and there is only one regular fixed point, so
we get a contradiction here as well.

\medskip

\item If $\rho(\pA_\ast) > 0$ then $\chi_{\pA_\ast}(X) = P_1(X)$. We
  have $\Tr(\pA_\ast)=-1$ and $\Tr(\pA^2_\ast)=3$, so that $L(\pA)=3$
  and $L(\pA^2)=-1$.  Since $L(\pA) = 3$ there are at least $3$ fixed
  singularities; thus we need only consider strata $(2,2,4)$ and
  $(2,2,2,2)$.

  $L(\pA)=3$ implies that all the singularities are necessarily fixed,
  with positive index. Let us denote by $\Sigma_1,\ \Sigma_2$ two
  degree-2 singularities.  Since $\textrm{Ind}(\pA,\Sigma_i)=1$, by
  Corollary~\ref{cor:trick} one has $\textrm{Ind}(\pA^2,\Sigma_i)=-3$,
  leading to $L(\pA^2) \leq -6 + 2 = -4$; but $L(\pA^2) = -1$, which
  is a contradiction.

\end{enumerate}

\subsection{Second polynomial: $\lambda(\pA) = \rho(P_2)$}

As in the previous section, we can rule out most strata associated
with~$P_2$ both for positive ($P_2(X)$) or negative ($P_2(-X)$)
dominant root.  For~$P_2(-X)$, however, there remain three strata that
cannot be eliminated:
\begin{equation*}
  (8),\ (2,6), \text{ and } (2,2,2,2).
\end{equation*}
We single out the last stratum, $(2,2,2,2)$, to illustrate that this
is a candidate.  Indeed, assume that three of the degree~$2$
singularities are cyclically permuted, and the fourth one is fixed.
For the triplet of singularities assume that the two ingoing (or
outgoing) separatrices are permuted by~$\phi^6$, so they are fixed
by~$\phi^{12}$.  At iterates~$3$ and~$9$ the three singularities are
fixed but their separatrices are permuted, and~$\rho(\phi^3)$
and~$\rho(\phi^9)$ are both negative, so by
Proposition~\ref{prop:index} the total contribution to the Lefschetz
number from these three singularities is 3.  At iterate~$6$ we
have~$\rho(\phi^6)>0$ but the separatrices are permuted, so again from
Proposition~\ref{prop:index} the total contribution is 3.  Finally, at
iterate~$12$ the singularities and their separatrices are fixed, so
the total contribution to~$L(\phi^{12})$ is~$3\cdot(1-4)=-9$.

For the fixed singularity of degree 2, assume that the two separatrices
are permuted by~$\phi^2$, so they are fixed by~$\phi^4$.  Hence, the
singularity contributes~$1$ to~$L(\phi^n)$ except when~$n$ is a
multiple of~$4$: we then have~$\rho(\phi^n)>0$ again by
Proposition~\ref{prop:index} the contribution is~$1-4=-3$.  As can be
seen in Table~\ref{tab:genus3-P2on2222}, the deficit in~$L(\phi^n)$
can be exactly made up by introducing regular periodic orbits (it is
easy to show that this can be done for arbitrary iterates).
\begin{table}
\begin{tabular}{c|ccccccccccccccc}
\hline
 $n$ & 1 & 2 & 3 & 4 & 5 & 6 & 7 & 8 & 9 & 10 & 11 & 12 & 13 & 14 & 15
 \\
\hline
 $L(\phi^n)$ & 2 & 0 & 5 & -4 & 7 & -3 & 16 & -12 & 23 & -25 & 46 & -55 & 80 & -112 & 160  \\
\hline
$L(2^3)$ & 0 & 0 & 3 & 0 & 0 & 3 & 0 & 0 & 3 & 0 & 0 & -9 & 0 & 0 & 3 \\
 $L(2^1)$ & 1 & 1 & 1 & -3 & 1 & 1 & 1 & -3 & 1 & 1 & 1 & -3 & 1 & 1 & 1 \\
 $L_{\text{ro}}$ & 1 & -1 & 1 & -1 & 6 & -7 & 15 & -9 & 19 & -26 & 45 & -43 & 79 & -113 & 156  \\
\hline
\end{tabular}
\medskip
\caption{
  For the first 15 iterates of~$\phi$, contribution to the Lefschetz
  numbers from the 
  various orbits, for the polynomial~$P_2(-X)$ from
  Table~\ref{tab:disc6} on stratum~$(2,2,2,2)$.  The first row specifies
  the iterate of~$\phi$; the second the total Lefschetz number; the
  third the contribution from the three permuted degree-2
  singularities; the fourth the contribution from the fixed degree-2
  singularity; the fifth the contribution from the regular (degree 0)
  orbits.  Note that~$L(2^3)$, $L(2^1)$, and~$L_{\mathrm{ro}}$ sum to~$L$.}
\label{tab:genus3-P2on2222}
\end{table}
To complete the proof of~\ref{thm:1} for~$g=3$, it remains to be shown
that such a homeomorphism can be constructed.

\subsection{Construction of a pseudo-Anosov homeomorphism by Dehn twists}
\label{sec:g3Dehn}

\begin{figure}
  \includegraphics[height=.15\textheight]{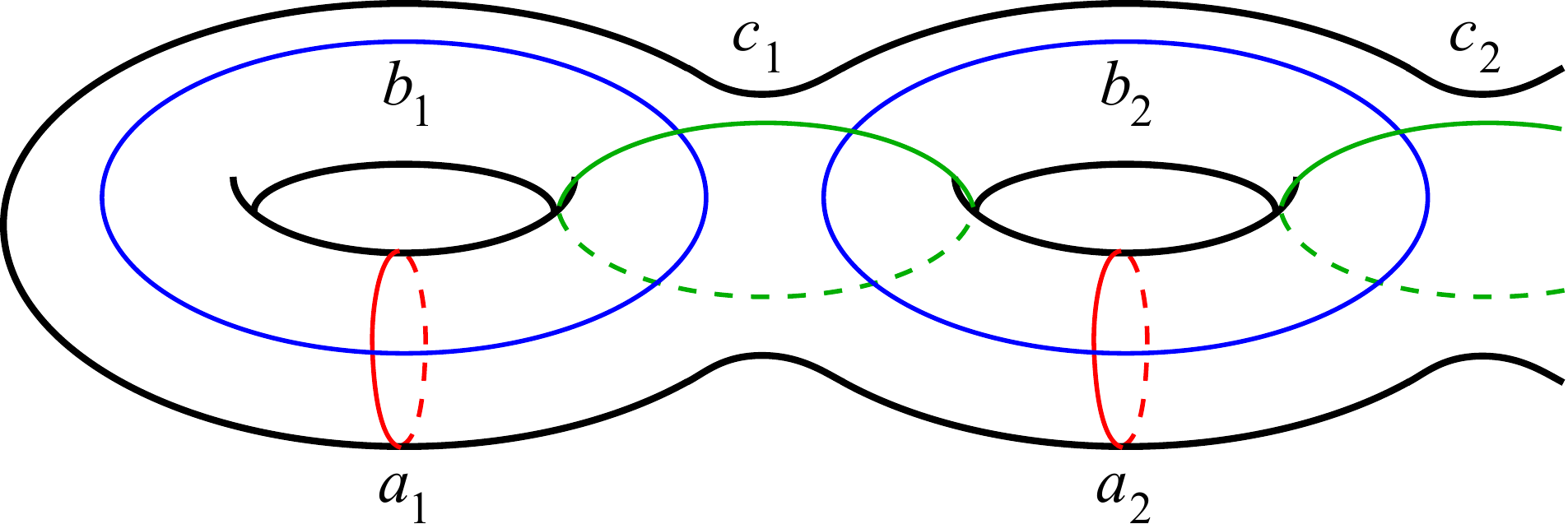}
\caption{Curves used to define Dehn twists.}
\label{fig:lickorish}
\end{figure}

We show how to realize in terms of Dehn twists a pseudo-Anosov
homeomorphism whose dilatation is the Perron root of~$P_2(X)$.  The
curves we use for Dehn twists are shown in
Figure~\ref{fig:lickorish}.
For example, a positive twist about~$c_1$ is written~$\DT{c_1}$; a
negative twist about~$b_2$ is written~$\DTi{b_2}$.  

\begin{Proposition}
There exists a pseudo-Anosov homeomorphism on a genus $3$ surface, stabilizing orientable foliations, and
having for dilatation the Perron root of the polynomial $P_2(X)$.
\end{Proposition}

\begin{proof}
Let us consider the sequence of Dehn twists
\begin{equation*}
  \DT{a_1}\DT{a_1}\DT{b_1}\DT{c_1}\DT{a_2}
  \DT{b_2}\DT{c_2}\DT{c_2}\DTi{a_3}\DTi{b_3}.
\end{equation*}
Its action on homology has~$P_2(-X)$ as a characteristic polynomial.
Since~$P_2(X)$ is irreducible and has no roots that are also roots of
unity~\cite{PisotSalem}, then by Theorem~\ref{thm:construction:pA} the
homeomorphism is isotopic to a pseudo-Anosov homeomorphism, say $f$
(we also use Bestvina's remark, Proposition~\ref{rem:Bestvina}).

We can compute the dilatation of $f$ by calculating the action on the
fundamental group (or using the code described in the remark below).
A straightforward calculation shows that the dilatation is the Perron
root of the polynomial $P_2(X)$, so~$f$ must also stabilize a pair of
orientable foliations.  Hence, it realizes our systole~$\delta_3^+$.
\end{proof}

\begin{Remark}
To search for pseudo-Anosov homeomorphisms, we used a computer code written by
Matthew D. Finn~\cite{FTJ2009}, which calculates the dilatation of
words in terms of Dehn twists. The code uses the fast method of
Moussafir~\cite{Moussafir2006} adapted to higher genus.  Hence, we can
examine a large number of words and find candidates with the required
dilatation.
\end{Remark}

\section{Genus two: A proof of Theorem~\ref{thm:genre2}}
\label{sec:genus2}

We prove theorem~\ref{thm:genre2} in two parts: we first find the
value of the systole~$\delta_2$, then demonstrate its uniqueness.

Recall that a surface $M$ of genus $g$ is called hyperelliptic if
there exists an involution $\tau$ (called the hyperelliptic
involution) with $2g+2$ fixed points. It is a classical fact that each
genus two surface is hyperelliptic. The fixed points are also called
the {\it Weierstrass points}. We now make more precise the
qualification ``up to hyperelliptic involution and covering
transformation'' of Theorem~\ref{thm:genre2}.

\begin{Remark}
\label{rem:unicity}
If $(M,q)$ is a hyperelliptic surface, then for each conjugacy class
of a pseudo Anosov homeomorphism $\pA$ on $M$ there exists another
conjugacy class, namely $\tau \circ \pA$, having the same dilatation.
For instance in genus $1$ the two Anosov homeomorphisms $\pA =
\left(\begin{smallmatrix}2&1\\1&1\end{smallmatrix}\right)$ and $\tau
\circ \pA =
\left(\begin{smallmatrix}-2&-1\\-1&-1\end{smallmatrix}\right)$ have
the same dilatation.

A second construction that produces another conjugacy class with the
same dilatation is the following. Let $\pA$ be a pseudo-Anosov
homeomorphism on a genus two surface $M$ stabilizing a non-orientable
foliation with singularity data $(1,1,2)$. Then there exists a
branched double covering $\pi : M \rightarrow \mathbb S^{2}$ such that
$\pA$ descends to a pseudo-Anosov $\tilde \pA$ on the sphere, fixing a
non-orientable measured foliation and having singularity data
$(-1,-1,-1,-1,-1,1,0)$ (see the proof of Theorem~\ref{thm:genre2}
below). Let the orientating double cover be $\pi' : N \rightarrow
\mathbb S^{2}$.  Now $\tilde \pA$ lifts to a new pseudo-Anosov
homeomorphism $\varphi$ on the genus-two surface $N$ (stabilizing
orientable foliations with singularity data $(4)$):
\[
\xymatrix{
    M \ar[r]^{\pA} \ar[d] & M \ar[d]  & N \ar[lld] \ar[r]^{\varphi} & N \ar[lld]\\
    \mathbb S^2 \ar[r]_{\tilde \pA} & \mathbb S^2
  }
\]
Now $\lambda(\pA) = \lambda(\tilde \pA) =\lambda(\varphi)$ (see
also~\cite{Lanneau} for more details). But of course the conjugacy
classes of $\pA$ and $\varphi$ are not the same.
\end{Remark}

Finally we will use the following result.

\begin{Proposition}
\label{prop:commute}
Let $(M,q)$ be a genus two flat surface and let $\tau$ be the affine
hyperelliptic involution. Let $\phi$ be an affine homeomorphism. Then
$\phi$ commutes with $\tau$.
\end{Proposition}

\begin{proof}[Proof of Proposition~\ref{prop:commute}]
  Let $\mathcal P =\{\fpq_1,\dots,\fpq_6\}$ be the set of Weierstrass
  points, i.e. the set of fixed points of $\tau$. \medskip

  Firstly let us show that $\pA$ preserves the set of Weierstrass
  points.  Since $\pA^{-1}\circ\tau\circ\pA$ is a non-trivial
  involution, it is an automorphism of the complex surface, thus the
  fixed points of $\pA^{-1}\circ\tau\circ\pA$ are also Weierstrass
  points. Let $p$ be a Weierstrass point. Then
  $\pA^{-1}\circ\tau\circ\pA(p)=p$ or $\tau\circ\pA(p)=\pA(p)$. Hence
  $\pA(p)$ is a fixed point of $\tau$, and thus $\pA(p)$ is a
  Weierstrass point. \medskip

  Now let $\psi=[\phi,\tau]=\phi \circ \tau\circ \phi^{-1}\circ\tau$
  be the commutator of $\phi$ and $\tau$. Since $\tau$ and $\phi$ are
  affine homeomorphisms, $\psi$ is also an affine homeomorphism. The
  derivative of $\psi$ is equal to the identity so that $\psi$ is a
  translation.  Since $\phi^{-1}\circ\tau(\fpq_1)=\phi^{-1}(\fpq_1)
  \in \mathcal P$ one has $\tau \circ\phi^{-1}\circ\tau(\fpq_1) =
  \phi^{-1}(\fpq_1)$ and
  $\psi(\fpq_1)=\phi\circ\phi^{-1}(\fpq_1)=\fpq_1$. The translation
  $\psi$ fixes a regular point. Thus it also fixes the separatrix
  issued from this point, and therefore $\psi=\id$ and $\phi$ commutes
  with $\tau$.
\end{proof}

\begin{proof}[Proof of Theorem~\ref{thm:genre2} (systole)]
  Let $\pA$ be a pseudo-Anosov homeomorphism with
  $\lambda(\pA)=\delta_2$. We know that $\delta^{+}_2$ is the Perron
  root of $X^4-X^3-X^2-X+1$ (see Zhirov~\cite{Zhirov1995}; see also
  Appendix~\ref{section:dehn} for a different construction).  Let us
  assume that $\delta_2 < \delta^{+}_2$. Thus $\pA$ preserves a pair
  of non-orientable measured foliations. The allowable singularity
  data for these foliations are $(2,2)$, $(1,1,2)$ or $(1,1,1,1)$.
  (Masur and Smillie~\cite{Masur:Smillie} showed that $(4)$ and
  $(1,3)$ cannot occur for non-orientable measured foliations.)
  \medskip

  It is well known that each genus two surface is a branched double
  covering of the standard sphere. Let $\pi : M \rightarrow \mathbb
  S^2$ be the covering and $\tau$ the associated involution.  It can
  be shown that $\tau$ is affine for the metric determined by $\pA$
  (see~\cite{Lanneau}).  Thus Proposition~\ref{prop:commute} applies
  and $\pA$ commutes with $\tau$. Hence $\pA$ induces a pseudo-Anosov
  homeomorphism $\tilde \pA$ on the sphere $\mathbb S^2$ with the same
  dilatation.  Of course $\tilde \pA$ leaves invariant a
  non-orientable pair of measured foliations.  The singularity data
  for $\pA$ are $(2,2)$, $(1,1,2)$, or $(1,1,1,1)$; The singularity
  data for $\tilde \pA$ are respectively $(-1,-1,-1,-1,0,0)$,
  $(-1,-1,-1,-1,-1,1,0)$, or $(-1,-1,-1,-1,-1,-1,1,1)$.  (For the
  first case, the singularity data cannot be $(-1,-1,-1,-1,-1,-1,2)$,
  otherwise the cover $\pi$ would be the orientating cover --- the
  branched points are precisely the singular points of odd degree, see
  Remark~\ref{rk:orientating:cover} --- thus the foliations of $\pA$
  would be orientable.)

  There exists an (orientating) double covering $\pi' : N \rightarrow
  \mathbb S^2$ such that $\tilde \pA$ lifts to a pseudo-Anosov
  homeomorphism $f$ on $N$ that stabilizes an orientable measured
  foliation. Actually, since the deck group is $\Z/ 2\Z$, there are
  two lifts: $f$ and $\tau \circ f$, where $\tau$ denote the
  hyperelliptic involution on $N$.  Since $\textrm{Tr}((\tau \circ
  f)_{\ast}) = - \textrm{Tr}(f)$, there is one lift, say $f$, with
  $\rho(\chi_{f_{\ast}}) > 0$.  By construction $\lambda(f) = \delta_2
  = \rho(\chi_{f_{\ast}})$.  Let us compute the genus of $N$ using the
  singularity data of $f$ as follows.

\begin{enumerate}

\item If the singularities of $\pA$ are $(2,2)$ then the singularities
  of $f$ are $(0)$; thus $N$ is a torus.

\item If the singularities of $\pA$ are $(1,1,2)$ then the
  singularities of $f$ are $(0,4)$; thus $N$ is a genus two surface.

\item If the singularities of $\pA$ are $(1,1,1,1)$ then the
  singularities of $f$ are $(4,4)$; thus $N$ is a genus three surface.

\end{enumerate}

In the first case one has $\delta_2 \geq \delta_1$, but since
$\delta_1>\delta^{+}_2$ this contradicts the assumption $\delta_2 <
\delta^{+}_2$.  In the second case $\delta_2 \geq \delta^{+}_2$ which
is also a contradiction. Let us analyze the third case.  Since
$\lambda(f) = \delta_2 < \delta^{+}_2$ and $f$ preserves an orientable
measured foliation on a genus three surface, Table~\ref{tab:genus2}
gives all possible minimal polynomials for $\delta_2$
with~$1<\rho(P)<\rho(X^4-X^3-X^2-X+1)$ (see
Appendix~\ref{sec:smallP}).
\begin{table}[htbp]
\begin{tabular}{lr}
\hline
polynomial  & Perron root  \\
\hline
$P_1 = (X^3-X-1)(X^3+X^2-1) $  & 1.32472 \\
$P_2 = X^6-X^4-X^3-X^2+1 $  & 1.40127\\ 
$P_3 = (X^3+X-1)(X^3-X^2-1) $  & 1.46557\\
$P_4 = X^6-X^5-X^3-X+1 $  & 1.50614\\
$P_5 = X^6-X^5-X^4+X^3-X^2-X+1 $  & 1.55603\\
$P_6 = X^6-2X^5+3X^4-5X^3+3X^2-2X+1 $  & 1.56769\\
$P_7 = X^6-X^4-2X^3-X^2+1 $  & 1.58235\\
$P_8 = X^6-2X^5+2X^4-3X^3+2X^2-2X+1$  & 1.63557\\
$P_9 = X^6-X^5+X^4-4X^3+X^2-X+1$  & 1.67114\\
\hline
\end{tabular}
\medskip
\caption{List of all reciprocal monic degree $6$ polynomials $P$ such 
that the Perron root $\lambda=\rho(P)$ satisfies $1 < \lambda <
\rho(X^4-X^3-X^2-X+1) \simeq 1.72208$.}
\label{tab:genus2}
\end{table}
We will obtain a contradiction for each case.  For each polynomial
$P_i$, we calculate the Lefschetz number of iterates of $f$ (see
Table~\ref{tab:lefchetz:genus:2}).
\begin{table}[htbp]
\begin{tabular}{llr}
\begin{tabular}{l|c c |}
\cline{2-3}
  & $L(f)$        & $L(f^2)$\\
\hline
\multicolumn{1}{|l|}{$P_1$} & $3$ &             \\
\multicolumn{1}{|l|}{$P_3$} &  & $3$             \\
\multicolumn{1}{|l|}{$P_6$}  & &    $4$          \\
\multicolumn{1}{|l|}{$P_9$}  & &    $3$        \\
\hline
\end{tabular} & \hskip 15mm &
\begin{tabular}{l|c c |}
\cline{2-3}
   & $L(f)$        &  $L(f^3)$        \\
\hline
\multicolumn{1}{|l|}{$P_2$} & $2$  &  $-1$           \\
\multicolumn{1}{|l|}{$P_4$} & $1$  &  $-2$           \\
\multicolumn{1}{|l|}{$P_5$}  & $1$         & $1$  \\
\multicolumn{1}{|l|}{$P_7$}  & $2$           & $-4$  \\
\hline
\end{tabular}
\end{tabular}
\medskip
\caption{Lefschetz number of iterates of
  the pseudo-Anosov homeomorphism~$f$.}
\label{tab:lefchetz:genus:2}
\end{table}
\begin{enumerate}

\item Polynomial $P_i$ for $i\in \{1,3,6,9\}$ cannot be a candidate 
since the number of singularities is $2$ and $L(f)$ or $L(f^2)$ is
greater than or equal to $3$.

\item Polynomial $P_i$ for $i\in \{2,4,5,7\}$ cannot be a
  candidate. Indeed the singularities are fixed with positive index,
  thus by Corollary~\ref{cor:trick} we should have $L(f^3) \leq -10$,
  but we know $L(f^3) \geq -4$ from Table~\ref{tab:lefchetz:genus:2}.

\end{enumerate}

Finally the last case we have to consider is $P_8$. In that case, the
Lefschetz number of $f$ is $0$ and the Lefschetz number of $f^3$ is
$-3$. Let $\Sigma_1$ and $\Sigma_2$ be the two singularities of $f$ on
$N$. Let us assume that the two singularities are fixed, so the index
of $f$ at $\Sigma_i$ is necessarily positive. Then by
Corollary~\ref{cor:trick} $\textrm{Ind}(f^3,\Sigma_i)=-5$, so that
$L(f^3)=-3=-10 -\#Fix(f^3)$ and $\#Fix(f^3)=-7$, which is a
contradiction. Hence $\Sigma_1$ and $\Sigma_2$ are exchanged by $f$,
and therefore by $f^3$. The formula $L(f^3)=-3$ reads $\#Fix(f^3)=3$,
so that $f$ has a unique length $3$ periodic orbit (and no fixed
points).  Recall also that $f$ commutes with the hyperelliptic
involution $\tau$ on $N$. This involution has exactly $8$ fixed points
on $N$: the two singularities and $6$ regular points, which we will
denote by $\{\Sigma_1,\Sigma_2,\fpq_1,\dots,\fpq_6\}$. \medskip

Let $\{\fpp,f(\fpp),f^2(\fpp)\}$ be the length-$3$ orbit. Since
$f\circ\tau=\tau \circ f$ the set
$\{\tau(\fpp),\tau(f(\fpp)),\tau(f^2(\fpp))\}$ is also a length-$3$
orbit and thus by uniqueness
\[
\left\{\fpp,f(\fpp),f^2(\fpp)\right\} =
\left\{\tau(\fpp),\tau(f(\fpp)),\tau(f^2(\fpp))\right\}.
\]
If $\tau(\fpp)=\fpp$ then $\fpp=\fpq_i$ for some $i$ and
$\{\fpp,f(\fpp),f^2(\fpp)\}$ is a subset of $\{\fpq_1,\dots,\fpq_6\}$.
Otherwise let us assume that $\tau(\fpp)=f(\fpp)$. Applying $f$ one
gets $f^2(\fpp)=f(\tau(\fpp))=\tau(f(\fpp))=\tau^2(\fpp)=\fpp$ which
is a contradiction. We get the same contradiction if
$\tau(\fpp)=f^2(\fpp)$. Therefore $\tau(\fpp)=\fpp$ and
$\{\fpp,f(\fpp),f^2(\fpp)\}$ is a subset of $\{\fpq_1,\dots,\fpq_6\}$.

Up to permutation one can assume that this set is
$\{\fpq_1,\fpq_2,\fpq_3\}$.  Since $f$ preserves the set
$\{\Sigma_1,\Sigma_2\}$ then $f$ also preserves
$\{\fpq_4,\fpq_5,\fpq_6\}$. Hence $f$ has a fixed point or another
length-$3$ periodic orbit, which is a contradiction. This ends the
proof of the first part of Theorem~\ref{thm:genre2}. \bigskip
\end{proof}

We now prove the uniqueness of the pseudo-Anosov homeomorphism
realizing the systole in genus two, up to conjugacy, hyperelliptic
involution, and covering transformations (see Remark~\ref{rem:unicity}).

\begin{proof}[Proof of Theorem~\ref{thm:genre2} (uniqueness)]
  We will prove that there is no other construction that realizes the
  systole in genus two.  The proof uses essentially McMullen's
  work~\cite{Mc}. Let $\pA$ and $\pA'$ be two pseudo-Anosov
  homeomorphisms on $M$ with $\lambda(\pA) = \delta_2$ and let
  $(M,q)$, $(M',q')$ be the two associated flat surfaces. \medskip

  The proof decomposes into $4$ steps. We first show that one can
  assume that $\pA$ and $\pA'$ leave invariant an orientable measured
  foliation with singularity data $(4)$. Then we show that we can
  assume, up to conjugacy, that the two surfaces $(M,q)$ and $(M',q')$
  are isometric. Finally we show that the derivatives $D \pA$ and $D
  \pA'$ of the affine homeomorphism on $M$ are conjugate. We then
  conclude that $\pA$ and $\pA$ are conjugated in the mapping class
  group Mod$(2)$.

  \begin{Step}
    If the foliation is non-orientable then we have seen (proof of
    Theorem~\ref{thm:genre2}) that the singularity data of $\pA$ is
    $(1,1,2)$. By Remark~\ref{rem:unicity} there exists a branched
    double covering $\pi : M \rightarrow \mathbb P^{1}$ such that
    $\pA$ descends to a pseudo-Anosov on the sphere $\mathbb P^{1}$
    with singularity $(-1,-1,-1,-1,-1,1,0)$. Now the orientating cover
    $\tilde \pi : \tilde M \rightarrow \mathbb P^{1}$ gives a
    pseudo-Anosov homeomorphism $\tilde \pA$ on the genus $2$ surface
    $\tilde M$, with orientable foliation and singularity data
    $(4)$. In addition $\lambda(\pA) = \lambda(\tilde \pA)$.  Hence,
    from this discussion one can assume that $\pA$ stabilizes an
    orientable measured foliation. The singularity data of the
    measured foliation is either $(4)$ or $(2,2)$. Using the Lefschetz
    theorem, one shows that $(2,2)$ is impossible.
  \end{Step}

  \begin{Step}
    Up to the hyperelliptic involution, we can assume that
    $\Tr(\pA)>0$ and $\Tr(\pA')>0$. There is natural invariant we can
    associate to a flat surface with a pseudo-Anosov homeomorphism
    $\varphi$: this is the trace field (see~\cite{Kenyon2000}), the
    number field generated by $\lambda(\varphi) +
    \frac1{\lambda(\varphi)}$. In our case of course the trace field
    of the surfaces $(M,q)$ and $(M',q')$ is the same since the
    dilatation of $\pA$ and $\pA'$ is the same. More precisely the
    trace field is $\mathbb Q[t]$, where $t = \delta_2 +
    \delta_2^{-1}$. A straightforward calculation gives that the
    minimal polynomial of $t$ is $X^2-X-3$, so the trace field is
    $\mathbb Q(\sqrt{13})$. \medskip

    Since the discriminant $\Delta=13\not \equiv 1$ mod $8$,
    Theorem~$1.1$ of~\cite{Mc} implies that there exists a $A\in
    \textrm{SL}_2(\mathbb R)$ such that $A(M,q)=(M',q')$. (We can
    always assume that the area of the flat surfaces $(M,q)$ and
    $(M',q')$ is $1$.) In particular there exists an affine
    homeomorphism $f : M \rightarrow M'$ such that $Df=A$. Hence
    $f^{-1} \pA' f$ is a pseudo-Anosov homeomorphism on the {\it same}
    affine surface $(M,q)$.
  \end{Step}

  \begin{Step}
    Now the derivatives of the two affine maps $\pA$ and $\pA'$ (on
    the same flat surface $(M,q)$) belong to the Veech group of the
    surface $(M,q)$. (This group has $3$ cusps and genus zero ---
    see~\cite{Mc}, Theorem~$9.8$.) Using the Rauzy--Veech induction,
    we can check that $D\pA$ and $A^{-1}D\pA' A$ are conjugated in
    this group.
  \end{Step}

  \begin{Step}
    Thus there exists $B\in \textrm{SL}_2(\mathbb R)$ such that
    $D\pA=B^{-1}D\pA' B$. Now let $h:M \rightarrow M$ be such that $D
    h=B$; hence one has $D\pA=Dh^{-1}D\pA'Dh$. Finally
    $h^{-1}\pA'h\pA^{-1}$ is an affine diffeomorphism with derivative
    map equal to the identity, and so it is a translation. Since the
    metric has a unique singularity (of type $(4)$),
    $h^{-1}\pA'h\pA^{-1}=\id$.  We conclude that $\pA$ and $\pA'$ are
    conjugate in the mapping class group Mod$(2)$, and the theorem is
    proved.
  \end{Step}
\end{proof}

\section{Genus four: A proof of Theorem~\ref{thm:1} for~$g=4$}

\subsection{Polynomials}

The techniques of the previous sections can also be applied to the
genus 4 case.  The only difference is that for genus four and higher
we rely on a set of Mathematica scripts to test whether a polynomial
is compatible with a given stratum.  This is straightforward: we
simply try all possible permutations of the singularities and
separatrices, and calculate the contribution to the Lefschetz numbers
for each iterate of~$\phi$.  Then we see whether the deficit in the
Lefschetz numbers can be exactly compensated by regular periodic
orbits.  If not, the polynomial cannot correspond to a pseudo-Anosov
homeomorphism on that stratum.

Again, we start with $\delta^{+}_4 \leq
\rho(X^8-X^7+X^6-X^5-X^4-X^3+X^2-X+1)\simeq 1.34372$ (for instance
see~\cite{Hironaka2006} or~\cite{Lanneau:Thiffeault}) and search for
candidate polynomials with smaller dilatation (see
Appendix~\ref{sec:smallP}), shown in
\begin{table}[htbp]
\begin{tabular}{lr}
\hline
polynomial & Perron root \\
\hline
$P_1 = X^8-X^5-X^4-X^3+1$ & 1.28064 \\
$P_2 = (X^3-X-1)(X^3+X^2-1)(X-1)^2$ & 1.32472 \\
$P_3 = (X^3-X-1)(X^3+X^2-1)(X+1)^2$ & 1.32472 \\
$P_4 = (X^3-X-1)(X^3+X^2-1)(X^2-X+1)$ & 1.32472 \\
$P_5 = (X^3-X-1)(X^3+X^2-1)(X^2+X+1)$ & 1.32472 \\
$P_6 = (X^3-X-1)(X^3+X^2-1)(X^2+1)$ & 1.32472 \\
\hline
\end{tabular}
\medskip
\caption{List of all reciprocal monic degree $8$ polynomials $P$ with 
Perron root $1 < \rho(P) < \rho(X^8-X^7+X^6-X^5-X^4-X^3+X^2-X+1)\simeq
1.34372$.
}
\label{tab:genus4}
\end{table}
Table~\ref{tab:genus4}.  Seeking a contradiction, we instead
immediately find that~$P_1(-X)$ is an allowable polynomial on strata
\begin{equation*}
  (2,10),\ (2,2,2,2,4), \text{ and } (2,2,2,6).
\end{equation*}
As an example we show the contributions to the Lefschetz numbers in
Table~\ref{tab:genus4-P1on10-2} on stratum~$(2,10)$.
\begin{table}
\begin{tabular}{c|ccccccccccccccc}
\hline
 $n$ & 1 & 2 & 3 & 4 & 5 & 6 & 7 & 8 & 9 & 10 & 11 & 12 & 13 & 14 & 15  \\
\hline
 $L(\phi^n)$ & 2 & 2 & 5 & -2 & 7 & -1 & 9 & -2 & 14 & -13 & 13 & -17 & 28 & -33 & 40 \\
\hline
 $L(10^1)$ & 1 & 1 & 1 & 1 & 1 & 1 & 1 & 1 & 1 & 1 & 1 & -11 & 1 & 1 & 1 \\
 $L(2^1)$ & 1 & 1 & 1 & -3 & 1 & 1 & 1 & -3 & 1 & 1 & 1 & -3 & 1 & 1 & 1 \\
 $L_{\text{ro}}$ & 0 & 0 & 3 & 0 & 5 & -3 & 7 & 0 & 12 & -15 & 11 & -3 & 26 & -35 & 38 \\
\hline
\end{tabular}
\medskip
\caption{For the first 15 iterates of~$\phi$, contribution to the
  Lefschetz numbers from the
  various orbits, for the polynomial~$P_1(-X)$ from
  Table~\ref{tab:genus4} on stratum~$(2,10)$.  See the caption to
  Table~\ref{tab:genus3-P2on2222} for details.}
\label{tab:genus4-P1on10-2}
\end{table}
Each singularity is fixed (as they must be since there is only one of
each type), and their separatrices are first fixed by~$\phi^{12}$
(degree~$10$) and~$\phi^4$ (degree~$2$), respectively.  We can easily
show that the Lefschetz numbers are consistent for arbitrary iterate.
It turns out that we can construct a pseudo-Anosov
homeomorphism having this dilatation.

\subsection{Construction of a pseudo-Anosov homeomorphism by Dehn twists}

We use the same approach as in Section~\ref{sec:g3Dehn} to find the
candidate word.

\begin{Proposition}
  There exists a pseudo-Anosov homeomorphism on a genus $4$ surface,
  stabilizing orientable foliations, and having for dilatation the
  Perron root of the polynomial $P_1(X)$.
\end{Proposition}

\begin{proof}[Proof]
Let us consider the sequence of Dehn twists
\begin{equation*}
  \DT{a_1}\DT{b_1}\DT{c_1}\DT{a_2}\DT{b_2}
  \DT{c_2}\DT{b_3}\DT{c_3}\DT{b_4}.
\end{equation*}
Its action on homology has~$P_1(-X)$ as a characteristic polynomial.
Since~$P_1(X)$ is irreducible and has no roots that are also roots of
unity~\cite{PisotSalem}, then by Theorem~\ref{thm:construction:pA} the
homeomorphism is isotopic to a pseudo-Anosov homeomorphism, say $f$.

We compute the dilatation of $f$ by calculating the action on the
fundamental group, which shows that the dilatation is the Perron root
of the polynomial $P_1(X)$.  Hence,~$f$ 
must also stabilize a pair of orientable foliations, and it
realizes our systole~$\delta_4^+$.
\end{proof}

\section{Higher genus}

\subsection{Genus five: A proof of Theorem~\ref{thm:1} for~$g=5$}
\label{sec:proofg5}

This time there is a known candidate with a lower dilatation than
Hironaka \& Kin's~\cite{Hironaka2006}: Leininger's pseudo-Anosov
homeomorphism~\cite{Leininger2004} having Lehmer's number~$\simeq
1.17628$ as a dilatation.  This pseudo-Anosov homeomorphism has
invariant foliations corresponding to stratum~$(16)$.  (The Lefschetz
numbers are also compatible with stratum~$(4,4,4,4)$.)  The polynomial
associated with its action on homology has~$\rho(P)<0$.  An exhaustive
search (see Appendix~\ref{sec:smallP}) leads us to conclude that there
is no allowable polynomial with a lower dilatation, so there is
nothing else to check.

As we finished this paper we learned that Aaber \&
Dunfield~\cite{Aaber2010} have found a pseudo-Anosov homeomorphism
with dilatation lower than~$\delta_5^+$ (stabilizing a non-orientable
foliation), implying that~$\delta_5 < \delta_5^+$.

\subsection{Genus six: A proof of Theorem~\ref{thm:lower} for~$g=6$
  (computer-assisted)}

For genus 6, we have demonstrated that the Lefschetz numbers
associated with~$P(-X)$, with~$P$ the polynomial in
Table~\ref{tab:thm:lower}, are compatible with stratum~$(16,4)$, with
Lehmer's number as a root (Lehmer's polynomial is a factor).  (There
is another polynomial with the same dilatation that is compatible with
the stratum~$(20)$.)  We have not yet constructed an explicit
pseudo-Anosov homeomorphism with this dilatation for genus 6, so
Theorem~\ref{thm:lower} is a weaker form than~\ref{thm:1}: it only
asserts that~$\delta^+_6$ is not less than this dilatation.  
Note, however, that whether or not this
pseudo-Anosov homeomorphism exists this is the first instance where
the minimum dilatation is not lower than for smaller genus.

\subsection{Genus seven: A proof of Theorem~\ref{thm:lower} for~$g=7$
  (computer-assisted)}
\label{sec:proofg7}

Again, we have not constructed the pseudo-Anosov homeomorphism
explicitly, but the Lefschetz numbers for the polynomial~$P(-X)$,
with~$P$ as in Table \ref{tab:thm:lower}, are compatible with
stratum~$(2,2,2,2,2,14)$.

As we finished this paper we learned that Aaber \&
Dunfield~\cite{Aaber2010} and Kin \& Takasawa~\cite{Kin2010} have
found a pseudo-Anosov homeomorphism with dilatation equal to the
systole~$\delta_7^+$.

\subsection{Genus eight: A proof of Theorem~\ref{thm:lower} for~$g=8$
  (computer-assisted)}

Genus eight is roughly the limit of this brute-force approach: it
takes our computer program about five days to ensure that we have the
minimizing polynomial.  The bound described in
Appendix~\ref{sec:smallP} yields~$5\times 10^{12}$ cases for the
traces, most of which do not correspond to integer-coefficient
polynomials.

Yet again, we have not constructed the pseudo-Anosov homeomorphism
explicitly, but the Lefschetz numbers for the polynomial~$P(-X)$,
with~$P$ as in Table~\ref{tab:thm:lower}, are compatible with
stratum~$(6,22)$.

As we finished this paper we learned that Hironaka~\cite{Hironaka2009}
has found a pseudo-Anosov homeomorphism with dilatation equal to the
systole~$\delta_8^+$.

\medskip


Examining the cases with even~$g$ leads to a natural question:

\begin{Question}
  Is the minimum value of the dilatation of pseudo-Anosov homeomorphisms
  on a genus $g$ surface, for $g$ even, with orientable invariant
  foliations, equal to the largest root of the
  polynomial~$X^{2g}-X^{g+1}-X^{g}-X^{g-1}+1$?
  \label{conj:systoleeven}
\end{Question}

\appendix

\section{Searching for polynomials with small Perron root}
\label{sec:smallP}

\subsection{Newton's formulas}

\mathnotation{\rbound}{\alpha}

The crucial task in our proofs is to find all reciprocal polynomials
with a largest real root bounded by a given value~$\rbound$ (typically
the candidate minimum dilatation).  Moreover, these must be allowable
polynomials for a pseudo-Anosov homeomorphism: the largest root (in
absolute value) must be real and strictly larger than all other roots,
and it must be outside the unit circle in the complex plane.

The simplest way to find all such polynomials is to bound the
coefficients directly.  For example, in genus 3, If we denote an
arbitrary reciprocal polynomial by
$P(X)=X^6+aX^5+bX^4+cX^3+bX^2+aX+1$, we want to find all polynomials
with Perron root smaller than $\rbound=\rho(X^3-X^2-1)\simeq 1.46557$
(the candidate minimum dilatation at the beginning of
Section~\ref{sec:proof:g3}).  Let $t=\rbound+\rbound^{-1}$; a
straightforward calculation assuming that half the roots of~$P(X)$ are
equal to~$\rbound$ shows
\[
  \lvert a\rvert \leq 3t, \qquad \lvert b\rvert\leq 3(t^2+1),\qquad
  \lvert c\rvert\leq t(t^2+6).
\]
Plugging in numbers, this means $\lvert a\rvert \le 6$, $\lvert
b\rvert \le 18$, and $\lvert c\rvert \le 26$.  Allowing for
$X\rightarrow -X$ since we only care about the absolute value of the
largest root, we have a total of~$12,765$ cases to examine.  Out of
these, only two polynomials actually have a root small enough and
satisfy the other constraints (reality, uniqueness of largest root),
as given in Section~\ref{sec:proof:g3}.

The problem with this straightforward approach (also employed by Cho
and Ham for genus 2, see~\cite{Ham2008}) is that it scales very poorly
with increasing genus.  For genus~$4$, the number of cases is
$9,889,930$; for genus~$5$, we have $63,523,102,800$ cases (we use
for~$\rbound$ the dilatation of Hironaka \& Kin's pseudo-Anosov
homeomorphism~\cite{Hironaka2006}, currently the best general upper
bound on~$\delta_g$).  As~$g$ increases, the target
dilatation~$\rbound$ decreases, which should limit the number of
cases, but the quantity~$t=\rbound+\rbound^{-1}$ converges to unity,
and the bound depends only weakly on~$\rbound-1$.

An improved approach is to start from Newton's formulas relating the
traces to the coefficients: for a
polynomial~$P(X)=X^{n}+a_{1}X^{n-1}+a_{2}X^{n-2}+\ldots +
a_{n-1}X+a_{n}$ which is the characteristic polynomial of a
matrix~$M$, we have
\begin{equation*}
  \Tr(M^k) =  \begin{cases}
    - k a_k -\sum_{m=1}^{k-1} a_m\, \Tr(M^{k-m}), \quad & 1 \le k \le n;\\
    - \sum_{m=1}^{n} a_m\, \Tr(M^{k-m}), \quad & k > n.
  \end{cases}
\end{equation*}
For a reciprocal polynomial, we have~$a_{n-k}=a_k$.  We can use these
formulas to solve for the~$a_k$ given the first few traces~$\Tr(M^k)$,
$1\le k\le g$ ($g=n/2$, $n$ is even in this paper).  We also have
\begin{Lemma}
\label{lem:Trbound}
  If the characteristic polynomial~$P(X)$ of a matrix~$M$ has a
  largest eigenvalue with absolute value~$r$, then
\begin{equation*}
  \l\lvert \Tr(M^k)\r\rvert \le n\,r^k;
\end{equation*}
Furthermore, if~$P(X)$ is reciprocal and of even degree, then
\begin{equation*}
  \l\lvert \Tr(M^k)\r\rvert \le
    \tfrac12 n(r^k + r^{-k}).
\end{equation*}
\end{Lemma}

\begin{proof}
Obviously,
\begin{equation*}
  \l\lvert \Tr(M^k)\r\rvert =  \l\lvert \sum_{m=1}^n s_m^k \r\rvert
  \le \sum_{m=1}^n \lvert s_m\rvert^k \le n\, r^k
\end{equation*}
where~$s_k$ are the eigenvalues of~$M$.  If the polynomial is
reciprocal and~$n$ is even, then 
\begin{equation*}
  \l\lvert \Tr(M^k)\r\rvert =
  \l\lvert \sum_{m=1}^{n/2} (s_m^k+s_m^{-k}) \r\rvert
  \le \tfrac12 n(r^k + r^{-k}).
\end{equation*}%
\end{proof}

We now have the following prescription for enumerating allowable
polynomials, given~$n$ and a largest root~$\rbound$:
\begin{enumerate}
  \item Use Lemma~\ref{lem:Trbound} to bound the traces $\Tr(M^k) \in
    \mathbb{Z}$, $k=1,\ldots,n/2$;
  \item\label{item:step2} For each possible set of~$n/2$ traces, solve for the
    coefficients of the polynomial;
  \item\label{item:step3} If these coefficients are not all integers,
    move on to the next possible set of traces;
  \item If the coefficients are integers, check if the polynomial is
    allowable: largest eigenvalue real and with absolute value less
    than~$\rbound$, outside the unit circle, and nondegenerate;
  \item Repeat step~\ref{item:step2} until we run out of possible
    values for the traces.
\end{enumerate}
Let's compare with the earlier numbers for~$g=5$:
assuming~$\Tr(M)\ge0$, we have $7,254,775$ cases to try, which is
already a factor of~$10^4$ fewer than with the coefficient bound.
Moreover, of these~$7,194,541$ lead to fractional coefficients, and so
are discarded in step~\ref{item:step3} above.  This only
leaves~$60,234$ cases, roughly a factor of~$10^6$ fewer than with the
coefficient bound.  Hence, with this simple approach we can tackle
polynomials up to degree~$16$ ($g=8$).  More refined approaches will
certainly allow higher degrees to be reached.

A final note on the numerical technique: we use Newton's iterative
method to check the dominant root of candidate polynomials.  A nice
feature of polynomials with a dominant real root is that their graph
is strictly convex upwards for~$x$ greater than the root (when that
root is positive, otherwise for~$x$ less than the root).  Hence,
Newton's method is guaranteed to converge rapidly and uniquely for
appropriate initial guess (typically, 5 iterates is enough for about 6
significant figures).  If the method does \emph{not} converge quickly,
then the polynomial is ruled out.

\subsection{Mahler measures}
\label{subsection:mahler}

Another approach is to use the Mahler measure of a polynomial. If $P$
is a degree $2g$ monic polynomial that admits a Perron root, say
$\rbound$, then the Mahler measure of $P$ satisfies $M(P) \le
\rbound^g$.  Thus to list all possible polynomials with a Perron root
less than a constant $\rbound$, we just have to list all possible
polynomials with a Mahler measure less than $\rbound^g$. Such lists
already exist in the literature (for example in~\cite{Boyd1980}).

\section{Rauzy--Veech induction and pseudo-Anosov
 homeomorphisms}
\label{appendix:construction:pA}

In this section we recall very briefly the basic construction of
pseudo-Anosov homeomorphisms using the Rauzy--Veech induction (for
details see~\cite{Veech1982}, \S 8, and
\cite{Rauzy,Marmi:Moussa:Yoccoz}).  We will use this to construct the
minimizing pseudo-Anosov homeomorphisms in genus~$3$ and~$4$.

\subsection{Interval exchange map}

Let $I \subset \mathbb R$ be an open interval and let us choose a
finite partition of $I$ into $d\geq 2$ open subintervals $\{I_j, \
j=1,\dots,d \}$.  An interval exchange map is a one-to-one map $T$
from $I$ to itself that permutes, by translation, the subintervals
$I_j$. It is easy to see that $T$ is precisely determined by a
permutation $\pi$ that encodes how the intervals are exchanged, and a
vector $\lambda=\{\lambda_j\}_{j=1,\dots,d}$ with positive entries
that encodes the lengths of the intervals.

\subsection{Suspension data}

A suspension datum for $T$ is a collection of vectors
$\{\zeta_j\}_{j=1,\dots,d}$ such that
\begin{enumerate}
\item $\forall j \in \{1,\dots,d\},\ Re(\zeta_j)=\lambda_j$;

\item $\forall k$, $1 \leq k \leq d-1,\ Im(\sum_{j=1}^k \zeta_j)>0$;

\item $\forall k$, $1 \leq k \leq d-1,\ Im(\sum_{j=1}^k \zeta_{\pi^{-1}(j)})<0$.

\end{enumerate}

To each suspension datum $\zeta$, we can associate a translation
surface $(M,q)=M(\pi,\zeta)$ in the following way. Consider the broken
line $L_0$ on $\mathbb{C}=\mathbb R^2$ defined by concatenation of the
vectors $\zeta_{j}$ (in this order) for $j=1,\dots,d$ with starting
point at the origin (see Figure~\ref{fig:surface}). Similarly, we
consider the broken line $L_1$ defined by concatenation of the vectors
$\zeta_{\pi^{-1}(j)}$ (in this order) for $j=1,\dots,d$ with starting
point at the origin.  If the lines $L_0$ and $L_1$ have no
intersections other than the endpoints, we can construct a translation
surface $S$ by identifying each side $\zeta_j$ on $L_0$ with the side
$\zeta_j$ on $L_1$ by a translation. The resulting surface is a
translation surface endowed with the form $\D z^2$.

Let $I \subset M$ be the horizontal interval defined by $I =
(0,\sum_{j=1}^d \lambda_j) \times \{0\}$.  Then the interval exchange
map $T$ is precisely the one defined by the first return map to $I$ of
the vertical flow on $M$.

\subsection{Rauzy--Veech induction}

The Rauzy--Veech induction $\mathcal R(T)$ of $T$ is defined as the
first return map of $T$ to a certain subinterval $J$ of $I$
(see~\cite{Rauzy,Marmi:Moussa:Yoccoz} for details).

We recall very briefly the construction.  The \emph{type}
$\varepsilon$ of $T$ is defined by $0$ if $\lambda_d >
\lambda_{\pi^{-1}(d)}$ and $1$ otherwise.  We define a subinterval $J$
of $I$ by
\[
J=\left\{ \begin{array}{ll}
I \backslash T(I_{\pi^{-1}(d)}) & \textrm{if $T$ is of type 0};\\
I \backslash I_d & \textrm{if $T$ is of type 1.} 
\end{array} \right.
\]
The Rauzy--Veech induction $\mathcal R(T)$ of $T$ is defined as the
first return map of $T$ to the subinterval $J$.  This is again an
interval exchange transformation, defined on $d$ letters (see
e.g.~\cite{Rauzy}). Moreover, we can compute the data of the new map
(permutation and length vector) by a combinatorial map and a matrix.
We can also define the Rauzy--Veech induction on the space of
suspensions.  For a permutation $\pi$, we call the \emph{Rauzy class}
the graph of all permutations that we can obtain by the Rauzy--Veech
induction. Each vertex of this graph corresponds to a permutation, and
from each permutation there are two edges labelled $0$ and $1$ (the
type). To each edge, one can associate a transition matrix that gives
the corresponding vector of lengths.

\subsection{Closed loops and pseudo-Anosov homeomorphisms}
\label{subsec:Veech}
We now recall a theorem of Veech:

\begin{NoNumberTheorem}[Veech]
  Let $\gamma$ be a closed loop, based at~$\pi$, in a Rauzy class and
  $R=R(\gamma)$ be the product of the associated transition
  matrices. Let us assume that $R$ is irreducible.  Let $\lambda$ be
  an eigenvector for the Perron eigenvalue $\alpha$ of $R$ and $\tau$
  be an eigenvector for the eigenvalue $\frac1{\alpha}$ of $R$. Then
\begin{enumerate}

\item $\zeta=(\lambda,\tau)$ is a suspension data for $T=(\pi,\lambda)$;

\item The matrix $A=\left( \begin{smallmatrix}\alpha^{-1} & 0 \\ 0 &
      \alpha \end{smallmatrix} \right)$ is the derivative map of an
  affine pseudo-Anosov diffeomorphism $\pA$ on the suspension
  $M(\pi,\zeta)$ over $(\pi,\lambda)$;

\item The dilatation of $\phi$ is $\alpha$;

\item All pseudo-Anosov homeomorphisms that fix a separatrix are
  constructed in this way.
\end{enumerate}

\end{NoNumberTheorem}

Since genus $4$ is simpler to construct than genus $3$, we present the
genus~$4$ case first in detail, and briefly outline the construction
of the other case.

\subsection{Construction of an example for $g=4$}

We shall prove

\begin{Theorem}
  There exists a pseudo-Anosov homeomorphism on a genus four surface,
  stabilizing orientable measured foliations, and having for
  dilatation the maximal real root of the polynomial
  $X^8-X^5-X^4-X^3+1$ (namely $1.28064...$).
\end{Theorem}

\subsubsection{Construction of the translation surface for $g=4$}

Let $|\alpha| > 1$ be the maximal real root of the polynomial
$P_1(X)=X^8-X^5-X^4-X^3+1$ with $\alpha<-1$, so that
$\alpha^8+\alpha^5-\alpha^4+\alpha^3+1=0$.  In the following, we will
present elements of $\mathbb Q[\alpha]$ in the basis
$\{\alpha^i\}_{i=0,\dots,7}$.  Thus the octuplet $(a_0,\dots,a_7)$
stands for~$\sum_{i=0}^7 a_i \alpha^i$.

\begin{figure}[htbp]
   \begin{center}
   \includegraphics[height=.45\textwidth]{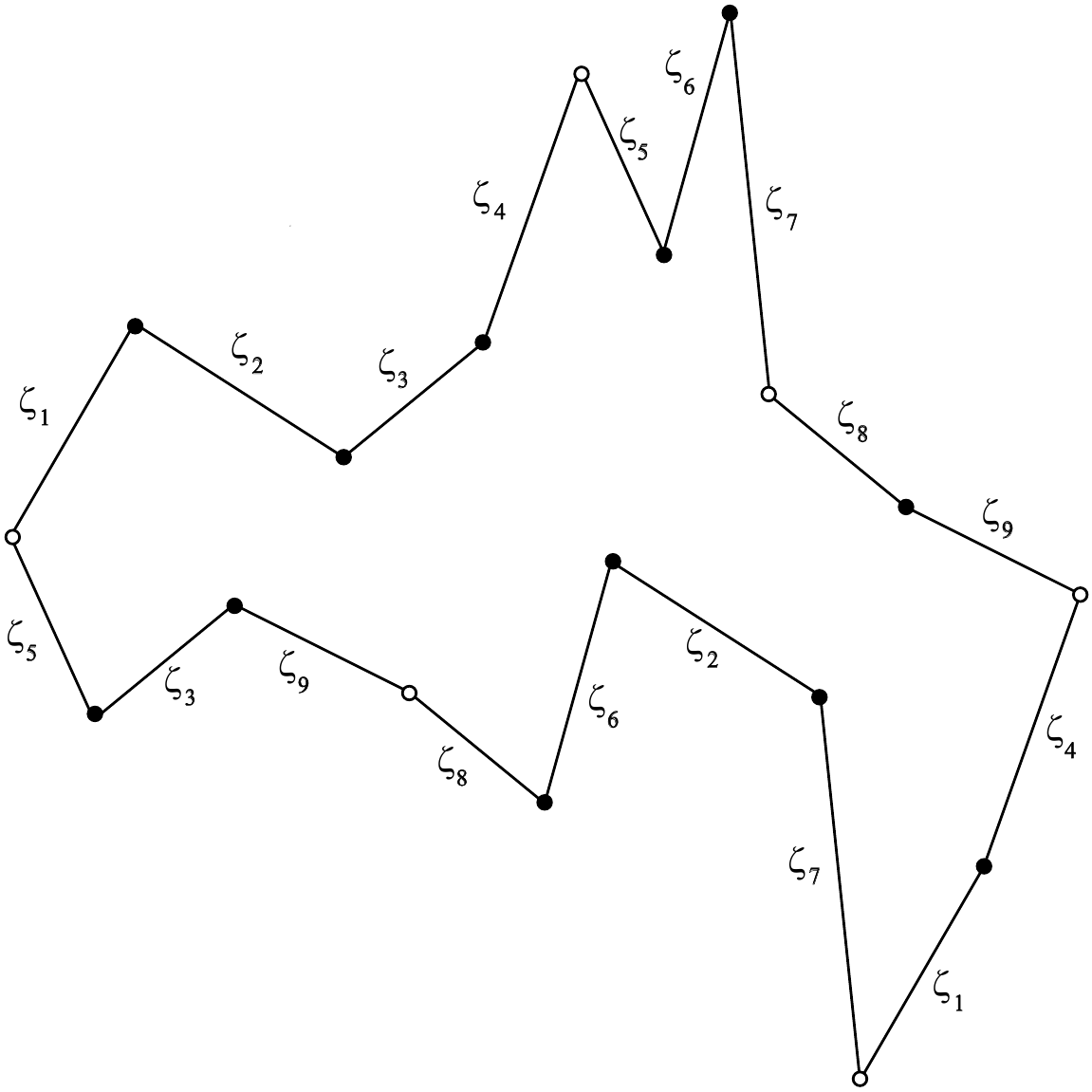}
   \caption{Construction of $(M,q)$. There are two singularities for
     the metric: one with conical angle $4\pi$ (hollow circles)
     and one with conical angle $12\pi$ (filled circles).  The
     stratum is thus $(2,10)$.}
     \label{fig:surface}
   \end{center}
\end{figure}

We start with the permutation $\pi=(5,3,9,8,6,2,7,1,4)$ and the closed
Rauzy path
\[
0-1-0-0-1-1-1-0-1-0-0-1-0-0.
\]
The associated Rauzy--Veech matrix is 
\[
R=\left(\begin{smallmatrix}
1&1&0&0&0&0&0&0&0 \\
0&0&1&1&1&1&1&1&0 \\
0&0&0&0&0&1&0&1&1 \\
1&0&0&1&0&0&1&0&0 \\
0&0&1&0&1&0&0&0&0 \\
0&0&0&1&1&0&0&0&0 \\
0&0&0&0&0&1&1&0&0 \\
1&0&0&0&0&1&1&1&0 \\
1&1&0&1&0&0&1&0&0 
\end{smallmatrix} \right).
\]
One checks that the characteristic polynomial of $R$ is $Q(X)$ with
the property that $Q(X)$ factors into $Q(X^4)=P_1(-X)S(X)$, where
$S(X)$ is a polynomial. Let $\lambda$ and $\tau$ be the corresponding
eigenvectors for the Perron root $\alpha^4$ of $Q$, expressed in the
$\alpha$-basis:
\[
\begin{array}{lcl}
\lambda_1 = (0,1,-2,1,-1,0,1,-1) \\
\lambda_2 = (0,-1,1,0,1,0,-1,0) \\
\lambda_3 = (-1,0,-1,0,0,-1,0,0) \\
\lambda_4 = (-1,2,-1,1,0,-1,1,0) \\
\lambda_5 = (1,-1,1,0,0,1,0,0) \\
\lambda_6 = (-1,1,-1,1,-1,-1,0,-1) \\
\lambda_7 = (1,-2,2,-2,1,1,-1,1) \\
\lambda_8 = (0,0,1,-1,1,0,0,1) \\
\lambda_9 = (1,0,0,0,0,0,0,0) 
\end{array}
\begin{array}{lcl}
\tau_1=(-1,0,0,0,0,-1,0,0)\\
\tau_2=(0,0,-1,1,0,1,0,-1)\\
\tau_3=(0,0,-1,0,-1,0,0,-1)\\
\tau_4=(0,1,0,0,0,0,1,0)\\
\tau_5=(0,0,0,1,0,0,0,0)\\
\tau_6=(0,0,-1,0,0,0,1,0)\\
\tau_7=(0,0,0,0,0,0,-1,0)\\
\tau_8=(0,1,0,0,0,0,0,0)\\
\tau_9=(-1,0,0,0,0,0,0,0).
\end{array}
\]
For $i=1,\dots,9$ we construct the vectors in $\R^{2}$ $\zeta_i
=\left( \begin{smallmatrix} \lambda_i \\ \tau_i \end{smallmatrix}
\right)$.  The resulting surface $(M,q)=M(\pi,\zeta)$ is drawn in
Figure~\ref{fig:surface}.

\subsubsection{Coordinates of the translation surface}

By construction, the coordinates of $(M,q)$ belong to $\mathbb
Q[\alpha]$.  We denote the vertices by $\fp_i$ for $i=1,\dots,18$ with
$\fp_1=0$ (see Figure~\ref{fig:surface:partition}).  Obviously for
$i\leq 9$, $\fp_i =\sum_{j=1}^i \zeta_j$, and for $i\geq 10$, $\fp_i
=\sum_{j=1}^9 \zeta_j - \sum_{j=1}^{i-9} \zeta_{\pi^{-1}(j)}$.  A
direct calculation gives
\[
\begin{array}{llll}
\fp_1 = ((0,0,0,0,0,0,0,0),& (0,0,0,0,0,0,0,0)) \\
\fp_2 = ((0, 1, -2, 1, -1, 0, 1, -1),& (-1, 0, 0, 0, 0, -1, 0, 0)  ) \\
\fp_3 = ((0, 0, -1, 1, 0, 0, 0, -1),&(-1, 0, -1, 1, 0, 0, 0, -1)   ) \\
\fp_4 = ( (-1, 0, -2, 1,0, -1, 0, -1) ,& (-1, 0, -2, 1, -1, 0, 0, -2)  )\\
\fp_5 = ((-2, 2, -3, 2, 0, -2, 1, -1) ,& (-1, 1, -2, 1, -1, 0, 1, -2)  )\\
\fp_6 = ((-1, 1, -2, 2, 0, -1, 1, -1) ,& (-1, 1, -2, 2, -1, 0, 1, -2)  )\\
\fp_7 = ((-2, 2, -3, 3, -1, -2, 1, -2) ,&  (-1, 1, -3, 2, -1, 0, 2, -2)  )\\
\fp_8 = ((-1, 0, -1, 1, 0, -1,0, -1) ,& (-1, 1, -3, 2, -1, 0, 1, -2)  )\\
\fp_9 = ((-1, 0, 0, 0, 1, -1, 0, 0) ,&  (-1, 2, -3, 2, -1, 0, 1, -2)  )\\
\fp_{10} =( (0, 0, 0, 0, 1, -1, 0, 0),&  (-2, 2, -3, 2, -1, 0,  1, -2) )\\
\fp_{11} = ((1, -2, 1, -1, 1, 0, -1, 0),& (-2, 1, -3, 2, -1, 0, 0, -2)  )\\
\fp_{12} = ((1, -3, 3, -2, 2, 0, -2, 1),&  (-1, 1, -3, 2, -1, 1, 0, -2) )\\
\fp_{13} =((0, -1, 1, 0, 1, -1, -1, 0) ,&  (-1, 1, -3, 2, -1, 1, 1, -2) )\\
\fp_{14} = ((0, 0, 0, 0, 0, -1, 0, 0),&  (-1, 1, -2, 1, -1, 0, 1, -1) )\\
\fp_{15} = ((1, -1, 1, -1, 1, 0, 0, 1),&  (-1, 1, -1, 1, -1,0, 0, -1) )\\
\fp_{16} = ((1, -1, 0, 0, 0, 0, 0, 0),&  (-1, 0, -1, 1, -1, 0, 0, -1)  )\\
\fp_{17} = ( (0, -1, 0, 0, 0, 0, 0, 0),& (0, 0, -1, 1, -1,  0, 0, -1)  )\\
\fp_{18} = ((1, -1, 1, 0, 0, 1, 0, 0),& (0, 0, 0, 1, 0, 0, 0, 0)  )
\end{array}
\]

\subsubsection{Construction of the pseudo-Anosov diffeomorphism}

Let $A$ be the hyperbolic matrix $\left( \begin{smallmatrix} \alpha^{-1} & 0\\
    0 & \alpha \end{smallmatrix} \right)$. Of course by construction
$A^4$ stabilizes the translation surface $(M,q)$ and hence there
exists a pseudo-Anosov homeomorphism on $M$ with dilatation
$\alpha^4$. We shall prove that this homeomorphism admits a root.

Let $(M',q')$ be the image of $(M,q)$ by the matrix $A$. We only need
to prove that $(M',q')$ and $(M,q)$ defines the same translation
surface, i.e. one can cut and glue $(M',q')$ in order to recover
$(M,q)$. This is

\begin{Theorem}
\label{theo:construc:genus4}
The surfaces $(M',q')$ and $(M,q)$ are isometric.
\end{Theorem}

\begin{Corollary}
  There exists a pseudo-Anosov diffeomorphism $f: X \rightarrow X$
  such that $D f=A$. In particular the dilatation of $f$ is
  $|\alpha|$.
\end{Corollary}

\begin{proof}[Proof of Theorem~\ref{theo:construc:genus4}]
Using the two relations $\alpha^8=-1-\alpha^3+\alpha^4-\alpha^5$ and
$\alpha^{-1}=\alpha^2-\alpha^3+\alpha^4+\alpha^7$ and the relations
that give the $\fp_i$, one gets by a straightforward calculation the
coordinates $\fq_i=A \fp_i$ of the surface $(M',q')$:
\[
\begin{array}{llll}
\fq_1 = ((0,0,0,0,0,0,0,0),& (0,0,0,0,0,0,0,0)) \\
\fq_{2} =((1, -2, 1, -1, 0, 1, -1, 0) ,& (0, -1, 0, 0, 0, 0, -1, 0)  )\\
\fq_{3} =( (0, -1, 1, 0, 0, 0, -1, 0) ,& (1, -1, 0, 0, 0, 1, 0, 0)  )\\
\fq_{4} =( (0, -2, 2, -1, 0,  0, -1, 1),&  (2, -1, 0, 0, -1, 1, 0, 0)  )\\
\fq_{5} =((2, -3, 4, -2, 0, 1, -1, 2) ,& (2, -1, 1, 0, -1, 1, 0, 1)  )\\
\fq_{6} =( (1, -2, 3, -1, 0, 1, -1, 1),&  (2, -1, 1, 0, 0, 1, 0, 1) )\\
\fq_{7} =( ( 2, -3, 5, -3, 0, 1, -2, 2),& (2, -1, 1, -1, 0, 1, 0, 2)  )\\
\fq_{8} =( (0, -1, 2, -1, 0, 0, -1, 1) ,& (2, -1, 1, -1, 0, 1, 0, 1)  )\\
\fq_{9} =( (0, 0, 1, 0, 0, 0, 0, 1) ,&  (2, -1, 2, -1, 0, 1, 0, 1) )\\
\fq_{10} =((0, 0, 0, 1, -1, 0, 0, 0) ,& (2, -2, 2, -1, 0, 1, 0, 1)  )\\
\fq_{11} =( (-2, 1, -2, 2, -1, -1, 0, -1),&  (2, -2, 1, -1, 0, 1, 0, 0) )\\
\fq_{12} =( (-3, 3, -3, 3, -1, -2, 1, -1),& (2, -1, 1, -1, 0, 1, 1, 0)  )\\
\fq_{13} =( (-1, 1, 0, 1, -1, -1, 0, 0) ,&(2, -1, 1, -1, 0, 1, 1, 1)   )\\
\fq_{14} =((0, 0, 0, 0, -1, 0, 0, 0) ,& (1, -1, 1, -1, 0, 0, 0, 1)  )\\
\fq_{15} =( (-1, 1, -2, 2, -1, 0, 1, -1),&  (1, -1, 1,  0, 0, 0, 0, 0)  )\\
\fq_{16} =((-1, 0, -1, 1, -1, 0, 0, -1) ,& (1, -1, 0, 0, 0, 0, 0, 0)  )\\
\fq_{17} =( (-1, 0, 0, 0, 0, 0, 0, 0) ,& (1, 0, 0, 0, 0, 0, 0, 0)  )\\
\fq_{18} =((-1, 1, -1, 1, 0, 0,0, -1) ,& (0,0, 0, 0, 1, 0, 0, 0)  )\\
\end{array}
\]

We will cut $M$ into several pieces in order to recover $M'$ such that
the boundary gluings agree.  Consider the decomposition in
Figure~\ref{fig:surface:partition}.
\begin{figure}[htbp]
   \begin{center}
     \includegraphics[height=10cm]{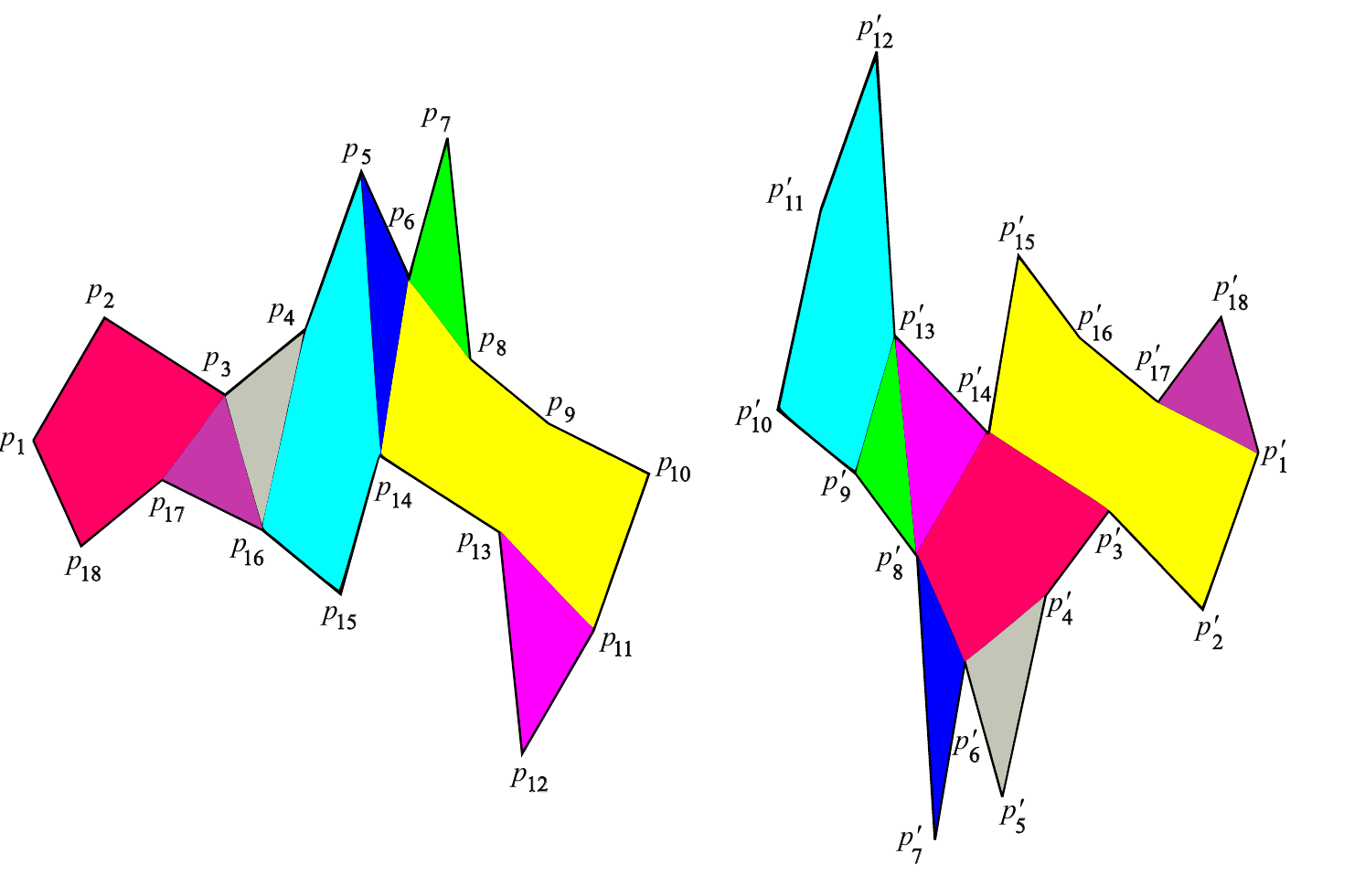}
     \caption{Partition of $(M,q)$ and
       $(M',q')=A(M,q)$.}
     \label{fig:surface:partition}
   \end{center}
\end{figure}
We enumerate the pieces on $M$ from the left to the right. For
instance, the first piece on $M$ has coordinates
$\fp_1\fp_2\fp_3\fp_{17}\fp_{18}$. The corresponding piece on $M'$ has
coordinates $\fq_3\fq_4\fq_6\fq_8\fq_{14}$. The translation is
$\overrightarrow{\fq_8\fp_1}= \overrightarrow{\fq_{14}\fp_2}=
\overrightarrow{\fq_3\fp_3}= \overrightarrow{\fq_4\fp_{18}}$.

\[
\begin{array}{clll}
\hline
\textrm{piece \#} & \textrm{coordinates on $M$} & \textrm{coordinates on $M'$} & \textrm{translation vectors}  \\
\hline
1&\fp_1\fp_2\fp_3\fp_{17}\fp_{18}&\fq_3\fq_4\fq_6\fq_8\fq_{14}&\overrightarrow{\fq_8\fp_1}= \overrightarrow{\fq_{14}\fp_2}=\overrightarrow{\fq_3\fp_3}= \overrightarrow{\fq_4\fp_{18}}\\[1pt]
2&\fp_3\fp_{16}\fp_{17}&\fq_{18}\fq_1\fq_{17}&\overrightarrow{\fq_{18}\fp_3}= \overrightarrow{\fq_{1}\fp_{16}}=\overrightarrow{\fq_{17}\fp_2}\\[1pt]
3&\fp_{3}\fp_{4}\fp_{16}&\fq_{6}\fq_{4}\fq_{5}&\overrightarrow{\fq_{6}\fp_{3}}= \overrightarrow{\fq_{4}\fp_{4}}=\overrightarrow{\fq_{5}\fp_{16}}\\[1pt]
4&\fp_{4}\fp_{5}\fp_{14}\fp_{15}\fp_{16}&\fq_{11}\fq_{12}\fq_{13}\fq_{9}\fq_{10}&\overrightarrow{\fq_{11}\fp_{4}}= \overrightarrow{\fq_{12}\fp_{5}}=
\dots=\overrightarrow{\fq_{10}\fp_{16}}\\[1pt]
5&\fp_{5}\fp_{6}\fp_{14}&\fq_{8}\fq_{6}\fq_{7}&\overrightarrow{\fq_{8}\fp_{5}}= \overrightarrow{\fq_{6}\fp_{6}}=\overrightarrow{\fq_{7}\fp_{14}}\\[1pt]
6&\fp_{6}\fp_{8}\fp_{9}\fp_{10}\fp_{11}\fp_{13}\fp_{14}&\fq_{15}\fq_{16}\fq_{17}\fq_{1}\fq_{2}\fq_{3}\fq_{14}&\overrightarrow{\fq_{15}\fp_{6}}= \overrightarrow{\fq_{16}\fp_{8}}=\dots=\overrightarrow{\fq_{14}\fp_{14}}\\[1pt]
7&\fp_6\fp_7\fp_8&\fq_8\fq_9\fq_{13}&\overrightarrow{\fp_6\fq_9}= \overrightarrow{\fp_7\fq_{13}}= \overrightarrow{ \fp_8\fq_8}\\[1pt]
8&\fp_{11}\fp_{12}\fp_{13}&\fq_{14}\fq_{8}\fq_{13}&\overrightarrow{\fq_{14}\fp_{11}}=
\overrightarrow{\fq_{8}\fp_{12}}=\overrightarrow{\fq_{13}\fp_{13}}\\
\hline
\end{array}
\]
\end{proof}

\subsection{Construction of an example for $g=3$} 

We shall prove

\begin{Theorem}
  There exists a pseudo-Anosov homeomorphism on a genus three surface,
  stabilizing orientable measured foliations, and having for
  dilatation the maximal real root of the polynomial $X^6 - X^4 - X^3
  - X^2 +1$ (namely $1.40127...$).
\end{Theorem}

\begin{proof}
  Let $|\alpha| > 1$ be the maximal real root of the polynomial
  $P_2(X)=X^6 - X^4 - X^3 - X^2 +1$ with $\alpha<-1$, so that
  $\alpha^6-\alpha^4+\alpha^3-\alpha^2+1=0$.  We start with the
  permutation $\pi=(6, 3, 8, 2, 7, 4, 10, 9, 5, 1)$ and the closed
  Rauzy path
\[
1-1-1-0-0-1-0-1-0-0.
\]
The associated Rauzy--Veech matrix is 
\[
R=\left(\begin{smallmatrix}
1&1&1&1&1&1&0&0&0&0 \\
0&0&0&0&0&0&1&0&0&0 \\
0&0&0&0&0&0&0&1&0&0 \\
0&0&0&0&0&0&0&0&1&0 \\
0&0&0&0&1&0&0&0&1&1 \\
0&0&1&0&0&1&0&0&0&0 \\
1&0&0&1&1&0&0&0&1&1 \\
0&1&0&0&0&0&0&0&0&0 \\
0&0&1&1&0&0&0&0&0&0 \\
0&0&0&0&1&1&0&0&0&0 
\end{smallmatrix} \right).
\]
The associated translation surface and its image are presented in
Figure~\ref{fig:surface:partition:genus3}.
\begin{figure}[htbp]
   \begin{center}
     \includegraphics[width=\textwidth]{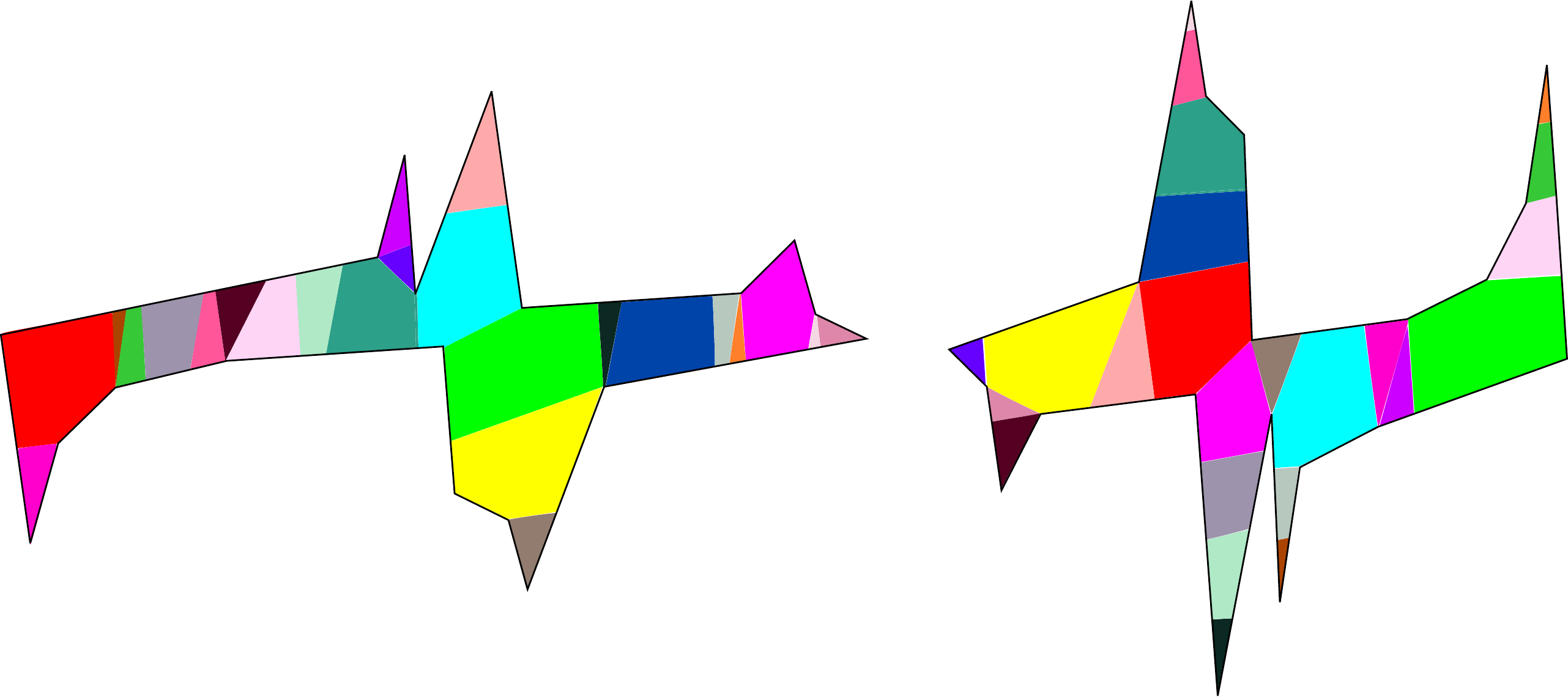}
     \caption{Partition of $(M,q)$ and $(M',q')=A(M,q)$. }
     \label{fig:surface:partition:genus3}
   \end{center}
\end{figure}
\end{proof}

\section{Genus two}
\label{section:dehn}

Let us consider the two sequences of Dehn twists on a genus two
surface,
\begin{equation*}
  T_{a_1}^2\DT{c_1}\DT{b_2}\DTi{a_2}\DT{b_1}
\qquad\text{ and }\qquad
  T_{a_1}^2\DTi{b_2}\DTi{c_1}\DTi{a_2}\DT{b_1}.
\end{equation*}

Their actions on the first homology group are respectively
$\left(\begin{smallmatrix} 1 & -3 & 0 & 1 \\ 1 & -2 & 0 & 1 \\ 0 & 2 & 2 & -1\\
    0 & 1 & 1 & 0 \end{smallmatrix} \right)$ and
$\left(\begin{smallmatrix} 1 & -1 & 1 & -1 \\ 1 & 0 & 1 & -1 \\ 0 & -1
    & -1 & 2\\
    0 & 0 & -1 & 1 \end{smallmatrix} \right)$. The characteristic
polynomials of these matrices are, respectively, $\xx^4 - \xx^3 -
\xx^2 - \xx +1$ and $\xx^4 - \xx^3 +3 \xx^2 - \xx +1$; thus
Theorem~\ref{thm:construction:pA} implies that the isotopy classes of
these homeomorphisms are pseudo-Anosov. Let $\pA_{1}$ and $\pA_{2}$ be
the corresponding maps.  One can calculate their dilatations from
their action on the fundamental group~\cite{FLP}. We check that the
dilatations, $\lambda$, are the same, namely the Perron root of the
polynomial $\xx^4 - \xx^3 - \xx^2 - \xx +1$ ($\lambda \simeq
1.72208$).

Theorem~\ref{theo:orientation:converse} thus implies that $\pA_{1}$
fixes an orientable measured foliation, and hence $\delta_{2}^{+} =
\lambda(\pA_{1})$ and $\pA_{2}$ fixes a non-orientable measured
foliation. We conclude that $\delta_{2}^{-} = \lambda(\pA_{2})$.

These two homeomorphisms are related by covering transformations (see
Remark~\ref{rem:unicity}).


\end{document}